\documentclass[review,3p,square]{elsarticle}%review

\usepackage{rotating}
\usepackage{amssymb}
\usepackage{amsmath}
\usepackage{amsthm}
\usepackage{stmaryrd}
\usepackage{bbm}
\usepackage{enumitem}

\usepackage{array}
\usepackage{multirow}
\usepackage{makecell}
\usepackage{longtable}
\usepackage{caption}

\usepackage{fancyhdr}
\usepackage{hyperref}
\hypersetup{%
  pdftitle={BSDEs},%
  pdfsubject={BSDEs},%
  pdfauthor={ShengJun FAN, Ying Hu},%
  pdfkeywords={BSDEs},%
  pdfstartview=FitH,%
  CJKbookmarks=true,%
  bookmarksnumbered=true,%
  bookmarksopen=true,%
  colorlinks=true, linkcolor=blue, urlcolor=blue, citecolor=blue, %
}

\usepackage{cleveref}
\crefname{thm}{Theorem}{Theorems}
\crefname{pro}{Proposition}{Propositions}
\crefname{lem}{Lemma}{Lemmas}
\crefname{rmk}{Remark}{Remarks}
\crefname{cor}{Corollary}{Corollaries}
\crefname{dfn}{Definition}{Definitions}
\crefname{ex}{Example}{Examples}
\crefname{section}{Section}{Sections}
\crefname{subsection}{Subsection}{Subsections}

\newcommand{\as}{{\rm d}\mathbb{P}\times{\rm d} t-a.e.}

\newcommand{\dif}{{\rm d}}
\newcommand{\ps}{\mathbb{P}-a.s.}

\newcommand{\essinf}{\mathop{\operatorname{ess\,inf}}}

\newcommand{\F}{\mathcal{F}}
\newcommand{\E}{\mathbb{E}}
\newcommand{\N}{\mathbb{N}}

\newcommand{\M}{\mathcal{M}}

\newcommand{\s}{\mathcal{S}}

\newcommand{\eps}{\varepsilon}
\newcommand{\T}{[0,T]}

\newcommand{\A}{\mathcal{A}}
\newcommand{\R}{{\mathbb R}}
\newcommand{\Q}{{\mathbb Q}}

\newtheorem{theorem}{Theorem}[section]
\newtheorem{lemma}[theorem]{Lemma}
\newtheorem{proposition}[theorem]{Proposition}
\newtheorem{remark}[theorem]{Remark}

\newtheorem{definition}[theorem]{Definition}
\newtheorem{example}[theorem]{Example}
\journal{arXiv}

\begin{document}
\begin{frontmatter}

\title{{On the uniqueness of solutions to quadratic BSDEs with non-convex generators and unbounded terminal conditions: the certain exponential moment case}\tnoteref{found}}
\tnotetext[found]{Partially supported by National Natural Science Foundation of China (No. 12171471).
\vspace{0.2cm}}

\author[Wang,Fan]{Yan Wang} \ead{wangyan\_shuxue@163.com}
\author[Fan]{Yaqi Zhang} \ead{TS22080028A31@cumt.edu.cn}
\author[Fan]{Shengjun Fan\corref{cor}} \ead{shengjunfan@cumt.edu.cn}\vspace{-0.6cm}

\affiliation[Wang]{organization={School of General Education, Jiangsu University of Architectural Technology},%Department and Organization
            %addressline={},
            city={Xuzhou 221116},
            %postcode={},
            %state={Jiangsu},
            country={China}}
            
\affiliation[Fan]{organization={School of Mathematics, China University of Mining and Technology},%Department and Organization
            %addressline={},
            city={Xuzhou 221116},
            %postcode={},
            %state={Jiangsu},
            country={China}}

\cortext[cor]{Corresponding \vspace{0.2cm}author}
\vspace{0.2cm}

\begin{abstract}
With the terminal value $|\xi|$ admitting some given exponential moments, we propose and prove several existence and uniqueness results for the unbounded solutions of quadratic backward stochastic differential equations whose generators may be represented as a uniformly continuous (not necessarily locally Lipschitz continuous) perturbation of some convex/concave function with quadratic growth. This perturbation satisfies various feasible conditions such as boundedness, sub-linear growth or linear growth. In particular, in some cases, the first component of the unique solution can be expressed as the value function of an optimal control problem. These results improves those posed in Delbaen, Hu and Richou [AIHP, 2011]
and Fan, Hu and Tang [2020, CRM] to some extent. The critical case is also tackled, which strengthens the main result of Delbaen, Hu and Richou [DCDS, 2011].\vspace{0.2cm}
\end{abstract}

\begin{keyword}
Backward stochastic differential equation \sep Existence and uniqueness \sep Unbounded solution\sep \\  \hspace*{1.85cm} Quadratic growth \sep Non-convex generator \vspace{0.2cm}

\MSC[2021] 60H10\vspace{0.2cm}
\end{keyword}

\end{frontmatter}
\vspace{-0.4cm}

\section{Introduction}
\label{sec:1-Introduction}
\setcounter{equation}{0}

In the whole paper, let us fix a positive integer $d$ and a positive real number $0<T<+\infty$. Let $(\Omega,\F,\mathbb{P})$ be a complete probability space carrying a standard $d$-dimensional Brownian motion $(B_t)_{t\geq0}, (\F_t)_{t\geq0}$ the completed natural $\sigma$-algebra generated by $(B_t)_{t\geq0}$ and $\F_T := \F$. We consider the following one-dimensional backward stochastic differential equation (BSDE in short):
\begin{align}\label{BSDE1.1}
     Y_t=\xi-\int_t^Tg(s,Z_s)\dif s+\int_t^TZ_s\cdot\dif B_s,\ \ \ t\in[0,T],
\end{align}
where $T$ is the terminal time, and the terminal value $\xi$ is an $\F_T$-measurable unbounded random variable, and the generator
\begin{align*}
g(\omega,t,z):\Omega\times[0,T]\times\R^d\mapsto\R
\end{align*}
is $(\F_t)$-progressively measurable for each $z\in\R^d$ which is continuous in $z$. BSDE with parameters $(\xi,T,g)$ is usually denoted by BSDE $(\xi,T,g)$. When $g$ has a quadratic growth in the state variable $z$, we call BSDE \eqref{BSDE1.1} a quadratic BSDE. By a solution of BSDE \eqref{BSDE1.1}, we mean a pair of $(\F_t)$-progressively measurable processes $(Y_t,Z_t)_{t\in[0,T]}$ valued in $\R\times\R^d$ such that $\ps$, the function $t\rightarrow Y_t$ is continuous, $t\rightarrow Z_t$ is square-integrable, $t\rightarrow g(t,Z_t)$ is integrable, and $(Y_\cdot,Z_\cdot)$ verifies \eqref{BSDE1.1}.\vspace{0.2cm}

Since the seminal paper \cite{Pardoux-Peng1990}, backward stochastic differential equations have become an active domain of research, and many applications have been found in various fields such as mathematical finance, partial differential equations (PDEs in short), stochastic control, nonlinear mathematical expectation and so on. In particular, a lot of efforts have been made in order to study the well-posedness of BSDEs, see for example \cite{Briand2003, Fan-Jiang2011, Fan-Jiang-Tian2011}. It should be noted that since it is widely used in the field of PDEs and financial mathematics, quadratic BSDEs have attracted much attention and are the subject of this paper. Interested readers are referred to for example \cite{Briand-Hu2006, Briand-Hu2008, Briand-Richou2017, Delbaen2011, Delbaen2015, Fan2016, Fan-Hu-Tang2020, Fan-Hu-Tang2023, Hu2005, Hu-Tang-Wang2022, Kobylanski2000, Richou2012, Rouge2000, Tian2023} for further details.\vspace{0.2cm}

We would like to especially mention the following works related closely to ours. Kobylanski \cite{Kobylanski2000} first studied the bounded solution of a quadratic BSDE with a bounded terminal value, and established a rather general existence and uniqueness result. The authors in \cite{Briand-Hu2006} obtained the first existence result for the unbounded solution of quadratic BSDEs with the terminal value $|\xi|$ admitting a certain exponential moment by applying the so-called localization method. Subsequently, uniqueness of the unbounded solution of quadratic BSDEs was established in \cite{Briand-Hu2008} when $|\xi|$ possesses all exponential moments and the generator $g$ is Lipschitz continuous in $y$, and either convex or concave with quadratic growth in $z$ (see (A1) and (A2) in Section \ref{Some existing existence and uniqueness results on quadratic BSDEs}). Based on the Legendre-Fenchel transform of convex functions, this result was further strengthened in \cite{Delbaen2011}. The main contribution of \cite{Delbaen2011} is to establish an existence and uniqueness result of the unbounded solution for a quadratic BSDE under the assumption that $|\xi|$ only admits a certain exponential moment. Furthermore, the critical case was addressed in \cite{Delbaen2015}, but an additional assumption that the generator $g$ is strongly convex in $z$ (see (A2') in Section \ref{Some existing existence and uniqueness results on quadratic BSDEs}) was required. We would like to especially mention that in all these articles mentioned above, uniqueness of the unbounded solution for a quadratic BSDE is obtained only when the generator $g$ is convex/concave in $z$. The case of a non-convex generator was tackled in \cite{Richou2012} and \cite{Briand-Richou2017}, but more assumptions are required on the terminal value $\xi$ than the exponential integrability. Recently, with $|\xi|$ possessing all exponential moments, the authors in \cite{Fan-Hu-Tang2020} proved a uniqueness result of the unbounded solution for a quadratic BSDE whose generator $g$ may be non-convex/non-concave in $z$, and instead satisfies a strictly quadratic condition (see (A1') in Section \ref{Some existing existence and uniqueness results on quadratic BSDEs}) and an extended convexity/concavity condition which holds typically for a locally Lipschitz perturbation of some convex/concave function (see (A3) in Section \ref{Some existing existence and uniqueness results on quadratic BSDEs}). Then, a question is naturally asked: when the terminal value $|\xi|$ only admits a certain exponential moment, and the quadratic growth generator $g$ may be represented as a locally Lipschitz perturbation of some convex/concave function, does the uniqueness for the unbounded solution of quadratic BSDEs hold still? The present paper gives some affirmative answers. Roughly speaking, it is verified in this paper that the uniqueness holds still even though the locally Lipschitz perturbation is replaced with a uniformly continuous (not necessarily locally Lipschitz continuous) one.\vspace{0.2cm}

More specifically, this paper is devoted to studying existence and uniqueness for the unbounded solution of quadratic BSDEs with non-convex/non-concave generators and some certain exponential moment conditions on the terminal value. We always assume that the generator $g(\omega,t,z):=g_1(\omega,t,z)+g_2(\omega,t,z)$, where $g_1$ is convex and has a quadratic growth in $z$ (see (A1) and (A2) in Section \ref{Some existing existence and uniqueness results on quadratic BSDEs}), and $g_2$ is only uniformly continuous (not necessarily locally Lipschitz continuous) in $z$ (see (H1) in Section \ref{Preliminaries}). With the terminal value $|\xi|$ admitting some given exponential moments, we propsoe and prove several existence and uniqueness results, which strengthen those established in \cite{Delbaen2011} and \cite{Fan-Hu-Tang2020} to some extent by relaxing the requirement of integrability on the terminal value and the convexity condition of the generator $g$ in $z$. Moreover, under an additional assumption that $g_1$ is strongly convex (see (A2') in Section \ref{Some existing existence and uniqueness results on quadratic BSDEs}), we also verify that uniqueness for the unbounded solution of quadratic BSDEs holds in the critical exponential moment case, which generalizes the main result obtained in \cite{Delbaen2015}.\vspace{0.2cm}

The remaining of this paper is organized as follows. In Section \ref{Preliminaries}, we introduce some basic notations, definitions, assumptions, lemmas, and propositions used later. Section \ref{Some existing existence and uniqueness results on quadratic BSDEs} collects some existing existence and uniqueness results on quadratic BSDEs posed in \cite{Delbaen2011,Delbaen2015,Fan-Hu-Tang2020}, providing a convenient reference for comparison with our results. The main results of this paper (Theorems \ref{ncth1} and \ref{ncth2}) are stated in Section \ref{Statement of the main results}. We work with the terminal value $|\xi|$ admitting an exponential moment of $p$-order for some $p>\gamma$ or just the critical case $p=\gamma$, where the constant $\gamma>0$ is defined in assumption (A1) of Section \ref{Some existing existence and uniqueness results on quadratic BSDEs}. With $g_1$ being convex in $z$, in (i) of Theorem \ref{ncth1} we prove uniqueness for the unbounded solution of BSDE \eqref{BSDE1.1} in the case of $g_2$ being bounded in $z$, while in (ii) of Theorem \ref{ncth1}, the boundedness of $g_2$ in $z$ is relaxed to have a sub-linear growth, but $g_1$ needs to additionally satisfy the strictly quadratic condition in $z$. And, in (iii) of Theorem \ref{ncth1}, $g_2$ only needs to  have a linear growth, but $g_1$ is supposed to be strongly convex in $z$. Then, in Theorem \ref{ncth2} we deal with the critical case where $g_1$ is strongly convex and $g_2$ is bounded in $z$. We compare our conclusions with related known results from \cite{Delbaen2011,Delbaen2015,Fan-Hu-Tang2020} in Remarks \ref{rmk:4.2} and \ref{rmk:4.4}. Furthermore, several examples to which Theorems \ref{ncth1} and \ref{ncth2} but no existing results can apply are provided in Example \ref{nclz}. We remark that our results also address quadratic BSDEs with non-concave generators (see (i) of Remark \ref{remark-the non-concave case} for more details). Finally, Sections \ref{Proof of Theorem 1} and \ref{Proof of Theorem 2} are devoted to the proofs of Theorems \ref{ncth1} and \ref{ncth2}, respectively. To do this, we borrow some ideas from \cite{Delbaen2011, Delbaen2015, Fan-Hu-Tang2020}, and systematically utilize some innovative techniques, including the Legendre-Fenchel transform of convex functions, Girsanov's theorem, the de La Vall$\acute{\rm{e}}$e Poussin lemma, Fenchel's inequality and Fatou's lemma. Furthermore, the existence and uniqueness on the $L^1$ solution of BSDEs along with the comparison theorem established in \cite{Fan2016} also plays an important role in our proof.\vspace{0.1cm}

\section{Preliminaries}\label{Preliminaries}
\label{sec:2-Preliminaries}
\setcounter{equation}{0}

In this section, we will first introduce some notations and spaces used in this paper, as well as the definition of $L^p~(p\geq1)$ solutions of BSDEs. Then, we will present an existence and uniqueness for the $L^1$ solution, a comparison theorem and two useful lemmas.\vspace{0.2cm}

The following are some notations that will be used later. Let $x\cdot y$ denote the Euclidean inner product for $x, y\in\R^d$. Denote by $\mathbbm{1}_A$ the indicator of set $A$, by $y^\top$ the transpose of vector $y$, and $\R_+:=[0,+\infty)$. Let $a\wedge b$ be the minimum of $a$ and $b$, $a^-:=-(a\wedge0)$ and $a^+:=(-a)^-$. For each $p>0$, $L^p(\Omega,\F_T,\mathbb{P})$ represents the set of $\F_T$-measurable random variables $|\xi|$ such that $\E[|\xi|^p]<+\infty$. We denote by $\partial f$ the subdifferential of a convex function $f:\R^{d}\rightarrow\R$, and the subdifferential of $f$ at $z_0\in\R^d$ is the non-empty convex compact set of elements $u\in\R^{1\times d}$ such that
$$
\forall z\in\R^d,\ \ \  f(z)-f(z_0)\geq u(z-z_0).
$$
Moreover, for any $(\F_t)$-progressively measurable $\R^{1\times d}$-valued process $(q_t)_{t\in[0,T]}$ such that $\ps$, $\int_0^T |q_s|^2 \dif s<+\infty$, we denote by $\mathcal{E}(q)$ the Dol{\rm{$\acute{e}$}}ans-Dade exponential
$$
\Big(\exp\Big(\int_0^tq_s\dif B_s-\frac{1}{2}\int_0^t |q_s|^2 \dif s\Big)\Big)_{t\in[0,T]}.
$$

Let us recall the following Fenchel's inequality
$$
xy\leq{\rm{exp}}(x)+y(\ln y-1),~~~\forall(x,y)\in\R\times(0,+\infty).
$$
Then for each $p>0$ we have
\begin{align*}
xy=px\frac{y}{p}\leq{\rm{exp}}(px)+\frac{y}{p}(\ln y-\ln p-1).
\end{align*}

Furthermore, we define the following two spaces of processes, where $p>0$.\vspace{0.2cm}

$\bullet$ $\s^p([0,T];\R)$ is the set of all $(\F_t)$-progressively measurable and continuous real-valued processes $(Y_t)_{t\in[0,T]}$ satisfying
$$\|Y\|_{{\s}^p}:=\left(\E\left[\sup_{t\in[0,T]} |Y_t|^p\right]\right)^{\frac{1}{p}\wedge1}<+\infty.$$
We set $\s=\bigcup_{p>1}\s^p$.\vspace{0.2cm}

$\bullet$ $\M^p([0,T];\R^d)$ is the set of all $(\F_t)$-progressively measurable $\R^d$-valued processes $(Z_t)_{t\in[0,T]}$ satisfying
$$
\|Z\|_{\M^p}:=\left\{\E\left[\left(\int_0^T |Z_t|^2\dif t\right)^{\frac{p}{2}}\right] \right\}^{\frac{1}{p}\wedge1}<+\infty.
$$
Recall that an $(\F_t)$-progressively measurable real-valued process $(Y_t)_{t\in[0,T]}$ belongs to class $(D)$ if the family of random variables $\{Y_\tau : \tau\in\Sigma_T\}$ is uniformly integrable where $\Sigma_T$ represents the set of all $(\F_t)$-stopping times $\tau$ valued in $[0,T]$. Throughout this paper, all equalities and inequalities between random variables are understood to hold $\ps$. For convenience, we recall the definition concerning the $L^p (p\geq1)$ solutions of BSDE \eqref{BSDE1.1}.\vspace{0.2cm}

\begin{definition}
Let generator $g(\omega,t,z):\Omega\times[0,T]\times\R^d\mapsto\R$ is $(\F_t)$-progressively measurable for each $z\in\R^d$, which is continuous in $z$. Assume that $(Y_\cdot,Z_\cdot)$ is a solution of BSDE \eqref{BSDE1.1}. If $(Y_\cdot,Z_\cdot)\in\s^p\times\M^p$ for some $p>1$, then it is called an $L^p$ solution; if $(Y_\cdot,Z_\cdot)\in\s^\beta\times\M^\beta$ for any $\beta\in(0,1)$ and $Y_\cdot$ belongs to class $(D)$, then an $L^1$ solution.\vspace{0.2cm}
\end{definition}

Now, we introduce an existence and uniqueness result for the $L^1$ solution of BSDE \eqref{BSDE1.1}, which is a direct consequence of Theorem 6.5 in \cite{Fan2016}.

\begin{proposition}\label{L1}
Assume that the terminal value $\xi+\int_0^T |g(t,0)|\dif t\in L^1(\Omega,\F_T,\mathbb{P})$, and the generator $g$ satisfies the following assumptions:

\begin{enumerate}
{\rm{\renewcommand{\theenumi}{\textbf{(H1)}}
\renewcommand{\labelenumi}{\theenumi}
\item $g$ is uniformly continuous in $z$, i.e., there exists a continuous nondecreasing function $\phi(\cdot):\R_+\rightarrow\R_+$ with $\phi(0)=0$ such that $\as$, for each $z_1, z_2\in{\R^{d}}$, we have
    $$
    |g(\omega,t,z_1)-g(\omega,t,z_2)|\leq \phi(|z_1-z_2|).
    $$
Without loss of generality, we can assume that $\phi(x)$ has a linear growth, i.e., there exists a constant $A>0$ such that $\phi(x)\leq A|x|+A$ for any $x\geq 0$.

\renewcommand{\theenumi}{\textbf{(H2)}}
\renewcommand{\labelenumi}{\theenumi}
\item $g$ has a sub-linear growth in $z$, i.e., there exist two nonnegative constants $a$ and $b$ along with a constant $\iota\in(0,1)$ such that $\as$, for each $z\in{\R^{d}}$, we have
$$ |g(\omega,t,z)|\leq a|z|^\iota+b.$$
}}
\end{enumerate}
Then BSDE \eqref{BSDE1.1} admits a unique $L^1$ solution.\vspace{0.2cm}
\end{proposition}

Furthermore, we present a comparison theorem that plays a key role later. This result corresponds to the known comparison results in Theorems 2.1 and 2.4 of \cite{Fan2016}, and we omit its proof here.\vspace{0.1cm}

\begin{proposition}\label{bj-nc}
Assume that $\xi$ and $\xi'$ are two terminal values, $g$ and $g'$ are two generators, and $(Y_\cdot,Z_\cdot)$ and $(Y'_\cdot,Z'_\cdot)$ are, respectively, a solution of BSDE $(\xi,T,g)$ and BSDE $(\xi',T,g')$.

\begin{enumerate}
\renewcommand{\theenumi}{{\rm{(i)}}}
\renewcommand{\labelenumi}{\theenumi}
\item If $g'$ satisfies assumption (H1), $(Y_\cdot-Y'_\cdot)^+\in\s$, $\ps$, $\xi\leq\xi'$ and
\begin{align}\label{9}
\as,\ \ \mathbbm{1}_{\{Y_t>Y'_t\}}\big(g(t,Z_t)-g'(t,Z_t)\big)\leq0,
\end{align}
then for each $t\in\T$, we have $\ps$, $Y_t\leq Y'_t$.

\renewcommand{\theenumi}{{\rm{(ii)}}}
\renewcommand{\labelenumi}{\theenumi}
\item Assume further that $\xi, ~\xi'\in L^1(\Omega, \F_T, P)$, and $(Y_\cdot,Z_\cdot)$ and $(Y'_\cdot,Z'_\cdot)$ are, respectively, an $L^1$ solution of BSDE $(\xi,T,g)$ and BSDE $(\xi',T,g')$. If $g'$ satisfies assumptions (H1) and (H2), $\ps$, $\xi\leq\xi'$, and \eqref{9} holds, then for each $t\in\T$, we have $\ps$, $Y_t\leq Y'_t$.\vspace{0.2cm}
\end{enumerate}
\end{proposition}

Let us close this section with following two lemmas. Lemma \ref{lem-poussin} is the de La Vall$\acute{\rm{e}}$e Poussin lemma, and Lemma \ref{lem-nc} comes from Lemma 3.2 in \cite{Delbaen2015}.\vspace{0.1cm}

\begin{lemma}\label{lem-poussin}
{\rm{(The de La Vall$\acute{\rm{e}}$e Poussin lemma):}} A family of random variables $\{X_i,i\in I\}$ is uniformly integrable if and only if there exists a nonnegative function $h:\R_+\rightarrow\R_+$ satisfying
$$
\lim_{x\rightarrow\infty}\frac{h(x)}{x}=+\infty\ \  {\rm and}\ \  \sup_{i\in I}\E[h(|X_i|)]<+\infty.
$$
In particular, if $\sup_{i\in I}\E[|X_i|\log |X_i|]<+\infty$, then $\{X_i,i\in I\}$ is uniformly integrable.\vspace{0.2cm}
\end{lemma}

\begin{lemma}\label{lem-nc}
Let $\mathcal{H}$ be a set of $(\F_t)$-progressively measurable random variables. The family of random variables $\{e^{\gamma X}|X\in\mathcal{H}\}$ is uniformly integrable if and only if there exists a strictly increasing differentiable function $k:\R_+\rightarrow\R_+$ such that $k(0)=\gamma, ~k(x)\rightarrow+\infty$ when $x\rightarrow+\infty$, and
$$
\sup_{X\in\mathcal{H}}\E[K(X^+)]<+\infty,
$$
with $K(x)=\int_0^x k(u)e^{\gamma u}\dif u, ~x\in\R_+$.\vspace{0.2cm}
\end{lemma}

\section{Some known results on quadratic BSDEs with unbounded terminal values}
\label{Some existing existence and uniqueness results on quadratic BSDEs}
\setcounter{equation}{0}

In this section, we collect some existing existence and uniqueness results for the unbounded solution of BSDE \eqref{BSDE1.1} reported in \cite{Delbaen2011,Delbaen2015,Fan-Hu-Tang2020}, and present two new observations. To do this, in this paper we are always given two constants $0<\bar {\gamma}\leq\gamma$, along with an $(\F_t)$-progressively measurable $\R_+$-valued process $(\alpha_t)_{t\in[0,T]}$. The following assumptions will be frequently used later, which come from \cite{Kobylanski2000, Briand-Hu2006, Briand-Hu2008, Delbaen2011, Fan-Hu-Tang2020, Delbaen2015}.\vspace{0.2cm}

\begin{enumerate}
\renewcommand{\theenumi}{\textbf{(A1)}}
\renewcommand{\labelenumi}{\theenumi}
\item $g$ has a quadratic growth in $z$, i.e., $\as$, for each $z\in\R^d$, we have
$$|g(\omega,t,z)|\leq \alpha_t(\omega)+\frac{\gamma}{2}|z|^2.$$

\renewcommand{\theenumi}{\textbf{(A2)}}
\renewcommand{\labelenumi}{\theenumi}
\item $g$ is convex in $z$, i.e., $\as$, $g(\omega, t, \cdot)$ is convex.\vspace{0.3cm}

\renewcommand{\theenumi}{\textbf{(A1')}}
\renewcommand{\labelenumi}{\theenumi}
\item $g$ satisfies a strictly-positive quadratic condition in $z$, i.e., $\as$, for each $z\in\R^d$, it holds that
$$
-\alpha_t(\omega)+\frac{\bar {\gamma}}{2}|z|^2\leq g(\omega,t,z)\leq\alpha_t(\omega)+\frac{\gamma}{2}|z|^2.\vspace{0.1cm}
$$

\renewcommand{\theenumi}{\textbf{(A2')}}
\renewcommand{\labelenumi}{\theenumi}
\item $g$ is strongly convex in $z$, i.e., $g$ satisfies assumption (A2), and there exist two constants $\varepsilon>0$ and $c\geq0$ such that $\as$, for each $z,z'\in\R^d, u_t\in\partial_z g(\omega,t,z')$,
\begin{align*}
g(\omega,t,z)-g(\omega,t,z')-u_t(\omega)(z-z')\geq\frac{\varepsilon}{2}|z-z'|^2-c.
\vspace{0.2cm}
\end{align*}
\end{enumerate}

\begin{remark}\label{rmk-js-nc}
{\rm{We have the following two remarks.

\begin{itemize}
\item [(i)] It is evident that (A1') implies (A1), and (A2') implies (A2).

\item[(ii)] It follows from Remark 1 in \cite{Delbaen2015} that if $\as, g(\omega,t,\cdot)$ is a second order differentiable function, then assumption (A2') is equivalent to the assumption: there exist $R\geq0$ and $\varepsilon>0$ such that $\as$, for all $z\in\R^d$ with $|z|>R$, we have $g''(\omega,t,z)\geq\varepsilon I_d$, where $I_d$ is a $d$-dimensional identity matrix.\vspace{0.2cm}
\end{itemize}}}
\end{remark}

The following proposition collects some results posed in \cite{Delbaen2011,Delbaen2015,Fan-Hu-Tang2020} on existence and uniqueness for the unbounded solution of a quadratic BSDE when the terminal value has a certain exponential moment. Note that the uniqueness in the following conclusion requires the generator $g$ to be convex, or even strongly convex.\vspace{0.1cm}

\begin{proposition}\label{cun}
Suppose that the terminal value $\xi$ satisfying
\begin{align}\label{xi-alpna}
\E\Big[\exp{\Big(p\Big(|\xi|+\int_0^T\alpha_t\dif t\Big)\Big)}\Big]<+\infty
\end{align}
for some $p\geq\gamma$, and the generator $g$ satisfies assumption (A1). We have

\begin{enumerate}
\renewcommand{\theenumi}{{\rm{(i)}}}
\renewcommand{\labelenumi}{\theenumi}
\item If $p>\gamma$, then BSDE \eqref{BSDE1.1} admits a solution $(Y_\cdot,Z_\cdot)$ such that
\begin{align}\label{kj}
\E\Big[\exp\Big(p\sup_{t\in[0,T]}\Big(|Y_t|+\int_0^t\alpha_s\dif s\Big)\Big)\Big]<+\infty,
\end{align}
and $Z_\cdot\in\M^2$. Morevoer, if $g$ further satisfies assumption (A2), then the solution $(Y_\cdot,Z_\cdot)$ of BSDE \eqref{BSDE1.1} satisfying \eqref{kj} is unique.\vspace{0.1cm}

\renewcommand{\theenumi}{{\rm{(ii)}}}
\renewcommand{\labelenumi}{\theenumi}
\item If $p=\gamma$, then BSDE \eqref{BSDE1.1} admits a solution $(Y_\cdot,Z_\cdot)$ such that the process $\big(e^{\gamma (|Y_t|+\int_0^t\alpha_s\dif s)}\big)_{t\in[0,T]}$ belongs to class $(D)$. Moreover, if $g$ also satisfies assumption (A2'), then the solution $(Y_\cdot,Z_\cdot)$ of BSDE \eqref{BSDE1.1} satisfying $\big(e^{\gamma (|Y_t|+\int_0^t\alpha_s\dif s)}\big)_{t\in[0,T]}$ belonging to class $(D)$ is unique.\vspace{0.1cm}
\end{enumerate}
\end{proposition}

\noindent\textbf{Proof.}
The result in (i) of Proposition \ref{cun} follows by exactly the same argument as in \cite{Delbaen2011}. The existence result in (ii) of Proposition \ref{cun} can be obtained by Proposition 1 of \cite{Fan-Hu-Tang2020}, and the uniqueness can be obtained by following the same idea in Theorem 4.1 of \cite{Delbaen2015} and applying its proof. The details are omitted \vspace{0.2cm}here. \hfill\framebox

The following Propositions \ref{Proposition 2 in [FanHuTang 2020]} comes from Proposition 2 in \cite{Fan-Hu-Tang2020}. It establishes an a priori estimate of quadratic variation of the first component in the unbounded solution of a BSDE with a strictly quadratic generator and an unbounded terminal value, and will be used later.

\begin{proposition}\label{Proposition 2 in [FanHuTang 2020]}
Let $\xi$ be a terminal condition, $g$ be a generator satisfying assumption (A1'), and $(Y_\cdot,Z_\cdot)$ be a solution to BSDE \eqref{BSDE1.1}. If
\begin{align*}
\E\Big[\exp\Big(p_0\sup_{t\in[0,T]}\Big(|Y_t|+\int_0^t\alpha_s\dif s\Big)\Big)\Big]<+\infty,
\end{align*}
for some real $p_0>0$, then there exists a constant $\eta>0$ depending only on $(\gamma, \bar {\gamma}, T, p_0)$ such that
\begin{align}\label{z1}
\E\Big[\exp\Big(\eta\int_0^T|Z_s|^2\dif s\Big)\Big]<+\infty.
\end{align}
In particular, for each $\lambda>0$, we have
\begin{align}\label{z2}
\E\Big[\exp\Big(\lambda\int_0^T|Z_s|\dif s\Big)\Big]<+\infty.
\end{align}
\end{proposition}
\vspace{0.2cm}

The following Proposition \ref{Theorem 5 in [FanHuTang 2020]} establishes two existence and uniqueness results for unbounded solutions of quadratic BSDEs where the terminal value admits exponential moments of any order, and the generator $g$ is not necessarily convex/concave in $z$. It is a direct consequence of Theorem 5 in \cite{Fan-Hu-Tang2020}, and we omit its proof here.

\begin{proposition}\label{Theorem 5 in [FanHuTang 2020]}
Suppose that the terminal value $\xi$ satisfies \eqref{xi-alpna} for any $p\geq\gamma$, and the generator $g$  satisfies assumption (A1) and the following assumption:\vspace{0.2cm}

\begin{enumerate}
\renewcommand{\theenumi}{\rm{\textbf{(A3)}}}
\renewcommand{\labelenumi}{\theenumi}
\item $g$ satisfies an extended convexity condition in $z$, i.e., there exist constants $\kappa\geq0$ and $\delta\in[0,1)$ such that $\as$, for each $z_1,z_2\in\R^d$ and $\theta\in(0,1)$,
$$
g(\omega,t,z_1)-\theta g(\omega,t,z_2)\leq (1-\theta)\Big(\alpha_t(\omega)+\kappa\big(|z_1|^{1+\delta}+|z_2|^{1+\delta}\big)+\gamma\Big|\frac{z_1-\theta z_2}{1-\theta}\Big|^2\Big).
$$
\end{enumerate}
Then, we have

\begin{enumerate}
\renewcommand{\theenumi}{{\rm{(i)}}}
\renewcommand{\labelenumi}{\theenumi}
\item If $\kappa=0$, then BSDE \eqref{BSDE1.1} admits a unique solution $(Y_\cdot,Z_\cdot)$ satisfying \eqref{kj} for any $p\geq\gamma$. Moreover, $Z_\cdot\in\M^p$ for any $p\geq1$.

\renewcommand{\theenumi}{{\rm{(ii)}}}
\renewcommand{\labelenumi}{\theenumi}
\item If $g$ also satisfies assumption (A1'), then BSDE \eqref{BSDE1.1} admits a unique solution $(Y_\cdot,Z_\cdot)$ satisfying \eqref{kj} for any $p\geq\gamma$. Moreover, \eqref{z1} holds for some real $\eta>0$.\vspace{0.2cm}
\end{enumerate}
\end{proposition}

\begin{remark}
{\rm According to Proposition 3 of \cite{Fan-Hu-Tang2020}, we know that a typical example of (A3) is
$$
g(\omega,t,z):=g_1(\omega,t,z)+g_2(\omega,t,z),
$$
where $g_1$ satisfies assumptions (A1) and (A2), and $g_2$ is $\delta$-locally Lipschitz continuous in $z$, i.e., there exists a constant $c\geq 0$ such that $\as$, for each $z_1,z_2\in\R^d$, we have
$$
|g_2(\omega,t,z_1)-g_2(\omega,t,z_2)|\leq c (1+|z_1|^\delta+|z_2|^\delta)|z_1-z_2|.
$$
In other words, a locally Lipschitz perturbation of some convex function satisfies assumption (A3). Consequently, Proposition \ref{Theorem 5 in [FanHuTang 2020]} strengthens the uniqueness result for quadratic BSDEs obtained in \cite{Briand-Hu2008}.}\vspace{0.2cm}
\end{remark}

The following Proposition \ref{observation 2} shows that assumption (A3) actually implies some kind of locally Lipschitz continuity of the generator $g$ in $z$.  This property is at the first time explored, to the best of our knowledge.\vspace{0.1cm}

\begin{proposition}\label{observation 2}
If the generator $g$ satisfies assumptions (A1) and (A3), then it is locally Lipschitz continuous in $z$, i.e., there exists a constant $C>0$ depending only on $(\gamma,\kappa,\delta)$ such that $\as$, for each $z_1,z_2\in \R^d$,
\begin{equation}\label{eq:3.9}
|g(\omega,t,z_1)-g(\omega,t,z_2)|\leq C(1+\alpha_t(\omega)+|z_1|^2+|z_2|^2)|z_1-z_2|.
\end{equation}
\end{proposition}

\noindent\textbf{Proof.} We only prove the case of $d=1$, based on which the general case can be easily verified. Now, suppose that the generator $g$ satisfies assumptions (A1) and (A3) with $d=1$. For each fixed $z_0\in \R$, let
$$
\tilde g(\omega,t,z):=g(\omega,t,z+z_0)-g(\omega,t,z_0),\ \ (\omega,t,z)\in \Omega\times\T\times\R.
$$
By (A3) and (A1), we have $\as$, $\tilde g(\omega,t,0)=0$ and for each $z_1,z_2\in \R$ and $\theta\in (0,1)$,
\begin{align}\label{eq:3.6}
\begin{split}
& \tilde g(\omega,t,(1-\theta)z_1+\theta z_2)-\theta \tilde g(\omega,t,z_2)\\
&\ \ =g(\omega,t,(1-\theta)(z_1+z_0)+\theta (z_2+z_0))-\theta g(\omega,t,z_2+z_0)+(1-\theta)g(\omega,t,z_0)\\
&\ \ \leq  (1-\theta)\left(\alpha_t(\omega)+\kappa (|z_1+z_0|^{1+\delta}+|z_2+z_0|^{1+\delta})+\gamma |z_1+z_0|^2\right)+(1-\theta)(\alpha_t(\omega)+\frac{\gamma}{2}|z_0|^2)\\
&\ \ \leq  (1-\theta)\bar C(1+\alpha_t(\omega)+|z_0|^2+|z_1|^2+|z_2|^2),
\end{split}
\end{align}
where $\bar C$ is a constant depending only on $(\gamma,\kappa,\delta)$. Letting $z_1=1, z_2=0$ and $z_1=-\theta, z_2=1-\theta$ in \eqref{eq:3.6} respectively yields that for each $\theta\in (0,1)$,
$$
\tilde g(\omega,t,1-\theta)\leq (1-\theta)\bar C(2+\alpha_t(\omega)+|z_0|^2)
$$
and
$$
-\theta\tilde g(\omega,t,1-\theta)\leq (1-\theta)\bar C(3+\alpha_t(\omega)+|z_0|^2),\vspace{0.2cm}
$$
which means that for each $\theta\in (1/2,1)$,
\begin{equation}\label{eq:3.7}
|\tilde g(\omega,t,1-\theta)|\leq 6\bar C (1+\alpha_t(\omega)+|z_0|^2)(1-\theta).
\end{equation}
Furthermore, letting $z_1=-1, z_2=0$ and $z_1=\theta, z_2=-(1-\theta)$ in \eqref{eq:3.6} yields that for each $\theta\in (0,1)$,
$$
\tilde g(\omega,t,-(1-\theta))\leq (1-\theta)\bar C(2+\alpha_t(\omega)+|z_0|^2)
$$
and
$$
-\theta\tilde g(\omega,t,-(1-\theta))\leq (1-\theta)\bar C(3+\alpha_t(\omega)+|z_0|^2),\vspace{0.2cm}
$$
which means that for each $\theta\in (1/2,1)$,
\begin{equation}\label{eq:3.8}
|\tilde g(\omega,t,-(1-\theta))|\leq 6\bar C (1+\alpha_t(\omega)+|z_0|^2)(1-\theta).
\end{equation}
Combining \eqref{eq:3.7} and \eqref{eq:3.8} yields that $\as$, for each $z\in (-1/2,1/2)$, we have
$$
|g(\omega,t,z+z_0)-g(\omega,t,z_0)|=|\tilde g(\omega,t,z)|\leq 6\bar C (1+\alpha_t(\omega)+|z_0|^2)|z|,
$$
for which the desired conclusion \eqref{eq:3.9} follows immediately for the case of $d=1$. \vspace{0.4cm}\hfill\framebox

The following proposition indicates that under certain mild conditions, assumption (A2') is strictly stronger than assumption (A1'), which will be used later.

\begin{proposition}\label{observation 1}
If there exists a process $u_t(\omega)\in\partial_z g(\omega,t,0)$ satisfying that
$$
\as,\ \ \ |u_t(\omega)|^2+|g(\omega,t,0)|\leq2\alpha_t(\omega),
$$
then assumptions (A1) and (A2') implies assumption (A1'), and the inverse implication is not true.\vspace{0.2cm}
\end{proposition}

\noindent\textbf{Proof.}
In fact, by letting $z'=0$ in (A2'), and combining Young's inequality, we have for each $t\in\T$,
\begin{align*}
g(\omega,t,z)&\geq\frac{\varepsilon}{2}|z|^2-c+g(\omega,t,0)+u_t(\omega)z\\
&\geq\frac{\varepsilon}{2}|z|^2-c-\Big(|g(\omega,t,0)|+\frac{1}{\varepsilon}|u_t(\omega)|^2\Big)-\frac{\varepsilon}{4}|z|^2\\
&\geq\frac{\varepsilon}{4}|z|^2-c-\Big(2+\frac{2}{\varepsilon}\Big)\alpha_t(\omega),
\end{align*}
which means that assumptions (A1) and (A2') implies assumption (A1'). Now, we can give an example to show that (A1') is strictly weaker than (A1) and (A2'). Define for each $z\in\R^d$,
\begin{align}\label{g-wan}
\hat{g}(z) := (2k-1)|z|-k(k-1),\ \ \ k-1\leq|z|<k,\ \ k\in\N.
\end{align}
It is not very difficult to check that for each $z\in\R^d$,
$$|z|^2\leq\hat{g}(z)\leq \frac{1}{4}+|z|^2,$$
and then $\hat{g}(z)$ satisfies assumption (A1'), but by virtue of the equivalent condition of (A2') mentioned in (ii) of Remark \ref{rmk-js-nc}, we can easily deduce that $\hat{g}(z)$ does not satisfy assumption (A2').\vspace{0.4cm}\hfill\framebox

We end this section with the following remark.

\begin{remark}\label{rmk:3.8}
{\rm It is evident that $(Y_\cdot,Z_\cdot)$ is a solution of BSDE $(\xi, T, g)$ if and only if $(-Y_\cdot,-Z_\cdot)$ is a solution of BSDE $(-\xi, T, \bar {g})$, where
$$
\bar {g}(\omega,t,z):=-g(\omega,t,-z),\ \ \ (\omega,t,z)\in \Omega\times\T\times\R^d.
$$
Moreover, when $\bar g$ respectively satisfies assumptions (A1), (A2), (A1'), (A2') and (A3), $g$ has a quadratic growth, is concave in $z$, satisfies a strictly-negative quadratic condition in $z$, is strongly concave in $z$, and satisfies an extended concavity condition in $z$. Thus, in light of preceding propositions and remarks reported in this section, some corresponding results on quadratic BSDEs with (strongly/extended) concave \vspace{0.2cm}generators can be obtained.}
\end{remark}

\section{Statement of the main results}
\label{Statement of the main results}
\setcounter{equation}{0}

This section will state several novel existence and uniqueness results on the unbounded solution of quadratic BSDEs. We consider the following one-dimensional BSDE:
\begin{align}\label{BSDE:4.1}
     Y_t=\xi-\int_t^Tg(s,Z_s)\dif s+\int_t^TZ_s\cdot\dif B_s,\ \ \ t\in [0,T].
\end{align}
We assume that the generator
\begin{align*}
g(\omega,t,z):\Omega\times[0,T]\times\R^d\mapsto\R
\end{align*}
is independent of the state variable $y$ and continuous in $z$, and that
$$
g(\omega,t,z):=g_1(\omega,t,z)+g_2(\omega,t,z),
$$
where $g_1$ is quadratic and convex in $z$ (i.e., (A1) and (A2)), and $g_2$ is uniformly continuous in $z$ (i.e., (H1)). Moreover, we also assume that for some $p\geq\gamma$,
\begin{align}\label{eq:4.2}
\E\Big[\exp{\Big(p\Big(|\xi|+\int_0^T\alpha_t\dif t\Big)\Big)}\Big]<+\infty.
\end{align}
In various cases, we further suppose that $g_1$ satisfies a strictly quadratic condition in $z$ (i.e., (A1')) or is strongly convex (i.e., (A2')), and that $g_2$ has a sub-linear growth in $z$ (i.e., (H2)) or isbounded (i.e., (H2) with $a=0$).\vspace{0.3cm}

The following Theorem \ref{ncth1} is the first main result of this paper, which establishes several existence and uniqueness results on the unbounded solution of quadratic BSDEs under the condition that the terminal value $\xi$ admits an exponential moment of $p$-order for some $p>\gamma$.\vspace{0.1cm}

\begin{theorem}\label{ncth1}
Suppose that $g:=g_1+g_2$ is a generator with $g_1$ satisfying assumptions (A1)-(A2) and $g_2$ satisfying assumption (H1), and that the terminal value $\xi$ satisfies \eqref{eq:4.2} for some $p>\gamma$. Then, 

\begin{enumerate}
\renewcommand{\theenumi}{{\rm{(i)}}}
\renewcommand{\labelenumi}{\theenumi}
\item If $g_2$ also satisfies assumption (H2) with $a=0$, then BSDE \eqref{BSDE:4.1}  admits a unique solution $(Y_\cdot,Z_\cdot)$ satisfying \eqref{kj}, and $Z_\cdot\in\M^2$. Moreover, we have $Y_\cdot=\essinf_{q\in\A}Y_\cdot^q$, and there exists $q^*\in\A$ such that $\as, ~Y_\cdot=Y_\cdot^{q^*}$, where $(Y_\cdot^q, Z_\cdot^q)$ is the unique $L^1$ solution of the following BSDE \eqref{bsde-qnc} under the probability measure $\Q_q$:
\begin{align}\label{bsde-qnc}
Y_t^q=\xi+\int_t^T\big(f_1(s,q_s)-g_2(s,Z_s^q)\big)\dif s+\int_t^T Z_s^q\cdot\dif B_s^q,~~~~t\in\T,
\end{align}
with $B_t^q:=B_t-\int_0^t q^\top_s\dif s, ~t\in[0,T]$ being a standard $d$-dimensional Brownian motion under $\Q_q$, the function $f_1$ being defined by
\begin{align}\label{fdy1}
f_1(\omega,t,q):=\sup_{z\in\R^d}\big(qz-g_1(\omega,t,z)\big),\ \  (\omega,t,q)\in\Omega\times\T\times\R^{1\times d},
\end{align}
and the admissible control set $\A$ being defined by
\begin{align*}
\A:=\bigg\{&(q_s)_{s\in[0,T]}~\text{is an } (\F_t)\text{-progressively measurable~}\R^{1\times d} \text{-valued process}:\\
&\int_{0}^{T}|q_s|^2\dif s<+\infty~ \ps,~~\E^{\Q_q}\Big[\int_{0}^{T}|q_s|^2\dif s\Big]<+\infty,\\
&\E^{\Q_q}\Big[|\xi|+\int_{0}^{T}|f_1(s,q_s)|\dif s\Big]<+\infty, ~\text{with} ~M_t^q:=\exp\Big(\int_{0}^{t}q_s\dif B_s-\frac{1}{2}\int_{0}^{t}|q_s|^2\dif s\Big), ~t\in[0,T] \\
&\text{ being a uniformly integrable martingale, and } \frac{\dif \Q_q}{\dif \mathbb{P}}:=M_{T}^q\bigg\}.\vspace{0.2cm}
\end{align*}

\renewcommand{\theenumi}{{\rm{(ii)}}}
\renewcommand{\labelenumi}{\theenumi}
\item If $g_1$ also satisfies assumption (A1'), and $g_2$ also satisfies assumption (H2), then BSDE \eqref{BSDE:4.1}  admits a unique solution $(Y_\cdot,Z_\cdot)$ satisfying \eqref{kj}, \eqref{z1} and \eqref{z2} for some $\eta>0$ depending only on $(\gamma, \bar {\gamma}, T, p)$ and each $\lambda>0$. Moreover, $Y_\cdot=\essinf_{q\in\A}Y_\cdot^q$, and there exists $q^*\in\A$ such that $\as, ~Y_\cdot=Y_\cdot^{q^*}$.

\renewcommand{\theenumi}{{\rm{(iii)}}}
\renewcommand{\labelenumi}{\theenumi}
\item If $g_1$ also satisfies assumption (A2'), and there exists a process $u_t(\omega)\in\partial_z g_1(\omega,t,0)$ such that $\as$,  $|u_t(\omega)|^2\leq\alpha_t(\omega)$, then BSDE \eqref{BSDE:4.1}  admits a unique solution $(Y_\cdot,Z_\cdot)$ satisfying \eqref{kj}, \eqref{z1} and \eqref{z2} for some $\eta>0$ depending only on $(\gamma, \bar {\gamma}, T, p)$ and each $\lambda>0$.\vspace{0.2cm}
\end{enumerate}
\end{theorem}

\begin{remark}\label{rmk:4.2}
{\rm With regard to Theorem \ref{ncth1}, we make the following remarks.

\begin{itemize}
\item[(i)] In (i) of Theorem \ref{ncth1}, $g_2$ is supposed to be bounded in $z$, while in (ii) of Theorem \ref{ncth1}, it is weakened to have a sub-linear growth, but $g_1$ needs to further satisfy a strictly quadratic condition in $z$. And, in (iii) of Theorem \ref{ncth1}, the requirement of $g_2$ in $z$ can be further weakened to have a linear growth, but $g_1$ needs to be strongly convex in $z$. Consequently, for the desired uniqueness and existence results, the growth condition of $g_2$ becomes weaker as the growth/convexity condition of $g_1$ becomes stronger.

\item[(ii)]  In (i) and (ii) of Theorem \ref{ncth1}, the first component of the unique solution of BSDE \eqref{BSDE:4.1} can be expressed as the value function of an optimal control problem, i.e., $Y_\cdot=\essinf_{q\in\A}Y_\cdot^q$. However, in (iii) of Theorem \ref{ncth1}, this conclusion dose not hold any longer since the $L^1$ solution $(Y_\cdot^q, Z_\cdot^q)$ of BSDE \eqref{bsde-qnc} under $\Q_q$ may not be well defined when $g_2$ only has a linear growth in $z$.

\item[(iii)]  In Theorem \ref{ncth1}, the terminal value $\xi$ is supposed to have an exponential moment of $p$-order for some $p>\gamma$. This condition is also required in (i) of Proposition \ref{cun}. However, it is supposed that the terminal value $\xi$ has exponential moments of any order in Proposition \ref{Theorem 5 in [FanHuTang 2020]}.

\item[(iv)] In Theorem \ref{ncth1}, the generator $g$ is allowed to be a convex function $g_1$ perturbed by a uniformly continuous function $g_2$, which is not necessarily locally Lipschitz continuous. However, it is supposed that in (i) of Proposition \ref{cun}, $g$ is convex in $z$, and that in Proposition \ref{Theorem 5 in [FanHuTang 2020]}, $g$ satisfies an extended convexity condition and then it must be locally Lipschitz continuous in $z$ according to Proposition \ref{observation 2}.

\item[(v)] It follows from (iii) and (iv) of Remark \ref{rmk:4.2} that Theorem \ref{ncth1} improves, to some extent, (i) of Proposition \ref{cun}, and Proposition \ref{Theorem 5 in [FanHuTang 2020]}.\vspace{0.1cm}
\end{itemize}}
\end{remark}

Theorem \ref{ncth1} shows existence and uniqueness of the unbounded solution for BSDE \eqref{BSDE:4.1} among solutions $(Y_\cdot, Z_\cdot)$ such that the inequality \eqref{kj} holds for some $p>\gamma$. The following theorem addresses the critical case: $p=\gamma$.

\begin{theorem}\label{ncth2}
Suppose that $g:=g_1+g_2$ is a generator with $g_1$ satisfying assumptions (A1) and (A2'), and $g_2$ satisfying assumptions (H1) and (H2) with $a=0$, and that the terminal value $\xi$ satisfies \eqref{eq:4.2} for $p=\gamma$, then BSDE \eqref{BSDE:4.1} admits a unique solution $(Y_\cdot,Z_\cdot)$ such that $(e^{\gamma(|Y_t|+\int_0^t\alpha_s\dif s)})_{t\in[0,T]}$ belongs to the class $(D)$.\vspace{0.2cm}
\end{theorem}

\begin{remark}\label{rmk:4.4}
{\rm We have the following comments.

\begin{itemize}
\item[(i)] In (ii) of Proposition \ref{cun}, the generator $g$ is supposed to be convex in $z$, while in Theorem \ref{ncth2}, $g$ is allowed to be a strongly convex function $g_1$ perturbed by a bounded uniformly continuous function $g_2$. Therefore, Theorem \ref{ncth2} strengthens (ii) of Proposition \ref{cun}. 

\item[(ii)] In Proposition \ref{Theorem 5 in [FanHuTang 2020]}, the terminal value $\xi$ is supposed to have exponential moments of any order, and the generator $g$ is supposed to be locally Lipschitz continuous in $z$ according to Proposition \ref{observation 2}, while in Theorem \ref{ncth2}, $\xi$ only requires the sharp exponential moment needed for the existence result, and $g$ may be not locally Lipschtiz continuous in $z$.

\item[(iii)] Compared with (iii) of Theorem \ref{ncth1}, in Theorem \ref{ncth2} the function $g_2$ is additionally supposed to be bounded, but the terminal value $\xi$ is relaxed to have a exponential moment of $\gamma$-order.\vspace{0.2cm}
\end{itemize}}
\end{remark}

\begin{example}\label{nclz}
{\rm We give several examples to which Theorems \ref{ncth1} or \ref{ncth2} but no existing result applies.

\begin{enumerate}
\renewcommand{\theenumi}{{\rm{(i)}}}
\renewcommand{\labelenumi}{\theenumi}
\item For each $(\omega,t,z)\in\Omega\times[0,T]\times\R^d$, define
$$
g(\omega,t,z):=|B_t(\omega)|+(1+\sin t)|z|^2+\mathbbm{1}_{0\leq|z|\leq1}\sqrt[4]{|z|}+\mathbbm{1}_{|z|>1}.
$$
It is easy to check that $g$ is non-convex with respect to $z$, but
$$g_1(\omega,t,z):=|B_t(\omega)|+(1+\sin t)|z|^2$$
satisfies (A1)-(A2), and
$$g_2(\omega,t,z):=\mathbbm{1}_{0\leq|z|\leq1}\sqrt[4]{|z|}+\mathbbm{1}_{|z|>1}\vspace{0.2cm}$$
satisfies (H1)-(H2) with
$$\alpha_\cdot\equiv|B_\cdot|, \ \gamma=4, \ A=1, \ a=0, \ b=1\ \  {\rm and} \ \ \phi(u)=\mathbbm{1}_{0\leq u\leq1}\sqrt[4]{u}+\mathbbm{1}_{u>1}.$$
It then follows from (i) of Theorem \ref{ncth1} that for each $\xi\in\F_T$ possessing an exponential moment of $p$-order with $p>4$, BSDE $(\xi,T,g)$ admits a unique solution $(Y_\cdot,Z_\cdot)$ satisfying \eqref{kj} and $Z_\cdot\in\M^2$.

\renewcommand{\theenumi}{{\rm{(ii)}}}
\renewcommand{\labelenumi}{\theenumi}
\item For each $(\omega,t,z)\in\Omega\times[0,T]\times\R^d$, define
$$
g(\omega,t,z):=\sqrt{|B_t(\omega)|}+\hat{g}(z)-\sqrt{|z|},
$$
where $\hat{g}(z)$ is defined in \eqref{g-wan}. Clearly,  $g$ is non-convex with respect to $z$, but $$g_1(\omega,t,z):=\sqrt{|B_t(\omega)|}+\hat{g}(z)$$
satisfies (A1') and (A2), and
$$g_2(\omega,t,z):=-\sqrt{|z|}\vspace{0.2cm}$$
satisfies (H1)-(H2) with
$$\alpha_\cdot\equiv\sqrt{|B_\cdot|}+1, \ \gamma=2, \ \bar {\gamma}=2, \ A=1, \ a=1, \ b=1, \ \iota=\frac{1}{2}\ \ {\rm and}\ \ \phi(u)=\sqrt{u}.$$
It then follows from (ii) of Theorem \ref{ncth1} that for each terminal value $\xi\in\F_T$ possessing an exponential moment of $p$-order with $p>2$, BSDE $(\xi,T,g)$ admits a unique solution $(Y_\cdot,Z_\cdot)$ satisfying \eqref{kj}, \eqref{z1} and \eqref{z2} for some $\eta>0$ and each $\lambda>0$.

\renewcommand{\theenumi}{{\rm{(iii)}}}
\renewcommand{\labelenumi}{\theenumi}
\item For each $(\omega,t,z)\in\Omega\times[0,T]\times\R^d$, define
$$
g(\omega,t,z):=\frac{1}{2}|z|^2+\mathbbm{1}_{0\leq|z|\leq1}\sqrt[3]{|z|^2}+\mathbbm{1}_{|z|>1}|z|.
$$
It is easy to check that $g$ is non-convex with respect to $z$, but
$$g_1(\omega,t,z):=\frac{1}{2}|z|^2$$
satisfies (A1) and (A2'), and $$g_2(\omega,t,z):=\mathbbm{1}_{0\leq|z|\leq1}\sqrt[3]{|z|^2}+\mathbbm{1}_{|z|>1}|z|$$
satisfies (H1) with
$$
\alpha_\cdot\equiv0,\  \gamma=1, \ \varepsilon=1, \ c=0,\  A=1\ \ {\rm and}\ \ \phi(u)=\mathbbm{1}_{0\leq u\leq1}\sqrt[3]{u^2}+\mathbbm{1}_{|u|>1}u.
$$
It then follows from (iii) of Theorem \ref{ncth1} that for each terminal value $\xi\in\F_T$ possessing an exponential moment of $p$-order with $p>1$, BSDE $(\xi,T,g)$ admits a unique solution $(Y_\cdot,Z_\cdot)$ satisfying \eqref{kj}, \eqref{z1} and \eqref{z2} for some $\eta>0$ and each $\lambda>0$.

\renewcommand{\theenumi}{{\rm{(iv)}}}
\renewcommand{\labelenumi}{\theenumi}
\item For each $(\omega,t,z)\in\Omega\times[0,T]\times\R^d$, define
$$
g(\omega,t,z):=\frac{1}{2}|z|^2-|z|+\mathbbm{1}_{|z|\leq\varepsilon} |z|\ln|z|+\mathbbm{1}_{|z|>\varepsilon}\varepsilon\ln\varepsilon
$$
with $\varepsilon\in(0,1)$ being small enough. It is easy to check that $g$ is non-convex in $z$, but $$g_1(\omega,t,z):=\frac{1}{2}|z|^2-|z|$$
satisfies (A1) and (A2'), and $$g_2(\omega,t,z):=\mathbbm{1}_{|z|\leq\varepsilon}|z|\ln|z|+
\mathbbm{1}_{|z|>\varepsilon}\varepsilon\ln\varepsilon$$
satisfies (H1)-(H2) with
$$
\alpha_\cdot\equiv1,\  \gamma=1,\  \varepsilon=1,\  A=2,\  a=0,\  b=2,\  c=1\ \  {\rm and}\ \  \phi(u)=\mathbbm{1}_{u\leq\varepsilon}u|\ln u|+\mathbbm{1}_{u>\varepsilon}\varepsilon|\ln\varepsilon|.
$$
It then follows from Theorem \ref{ncth2} that for each $\xi\in\F_T$ possessing an exponential moment of $1$-order, BSDE $(\xi,T,g)$ admits a unique solution $(Y_\cdot,Z_\cdot)$ such that $(e^{|Y_t|})_{t\in[0,T]}$ belongs to class $(D)$.
\end{enumerate}
}
\end{example}

We would like to especially emphsis that all $g_2$ in these examples are not locally Lipschitz continuous in $z$. It follows from Proposition \ref{observation 2} that Proposition \ref{Theorem 5 in [FanHuTang 2020]} cannot be applied to the prvious examples.\vspace{0.3cm}

For readers' convenience, we summarize the main results of Theorems \ref{ncth1} and \ref{ncth2} in the following table: Main results in Theorems \ref{ncth1} and \ref{ncth2}.

\begin{center}
  \captionof{table}{Main results in Theorems \ref{ncth1} and \ref{ncth2}}
  \vspace{-6pt}
  \begin{tabular}{|>{\centering\arraybackslash}m{2.0cm}|
                   >{\centering\arraybackslash}m{2.5cm}|
                   >{\centering\arraybackslash}m{3.4cm}|
                   >{\centering\arraybackslash}m{4.8cm}|
                   >{\centering\arraybackslash}m{2.0cm}|}
    \hline
    \multicolumn{1}{|>{\centering\arraybackslash}m{2.0cm}|}{Cond. on $|\xi|+\int_0^T \alpha_t\dif t$}
    & \multicolumn{1}{>{\centering\arraybackslash}m{2.5cm}|}{Cond. on $g_1$}
    & \multicolumn{1}{>{\centering\arraybackslash}m{3.4cm}|}{Cond. on $g_2$}
    & \multicolumn{1}{>{\centering\arraybackslash}m{4.8cm}|}{Space of solution}
    & \multicolumn{1}{>{\centering\arraybackslash}m{2.0cm}|}{Where} \\
    \hline
    \hspace*{-0.2cm}\multirow{8}{2.5cm}{\centering exponential moment of $p$-order for some $p>\gamma$}
    & quadratic and convex
    & uniformly continuous and bounded
    & $\sup\limits_{t\in[0,T]}\big(|Y_t|+\int_0^t\alpha_s\dif s\big)$ has an exponential moment of $p$-order and $Z_\cdot\in\M^2$
    & (i) of Theorem \ref{ncth1} \\
    \cline{2-5}
    & \vspace{0.5cm}strict quadratic and convex \vspace{0.5cm}
    & uniformly continuous and sub-linear
    & \multirow{2}{4.8cm}{\centering $\sup\limits_{t\in[0,T]}\big(|Y_t|+\int_0^t\alpha_s\dif s\big)$ has an exponential moment of $p$-order and $\E\Big[\exp\Big(\eta\int_0^T|Z_s|^2\dif s\Big)\Big]<+\infty$ for some $\eta>0$}
    & (ii) of Theorem \ref{ncth1} \\
    \cline{2-3} \cline{5-5}
    & \vspace{0.5cm}quadratic and strongly convex \vspace{0.5cm}
    & uniformly continuous
    & & (iii) of Theorem \ref{ncth1} \\
    \hline
    exponential moment of $\gamma$-order
    & quadratic and strongly convex
    & uniformly continuous and bounded
    & $(e^{\gamma(|Y_t|+\int_0^t\alpha_s\dif s)})_{t\in[0,T]}$ is of class $(D)$
    & Theorem \ref{ncth2} \\
    \hline
  \end{tabular}
\end{center}

\vspace{0.2cm}

\begin{remark}\label{remark-the non-concave case}
{\rm We end this section by the following two comments.
\begin{itemize}
\item[(i)] In light of Remark \ref{rmk:3.8}, Theorems \ref{ncth1} and \ref{ncth2} also address the unbounded solutions of quadratic BSDEs with generator $g$ being allowed to be a concave function $g_1$ perturbed by a uniformly continuous function $g_2$, which is not necessarily locally Lipschitz continuous in $z$. 

\item[(ii)] Note that when $g$ also depends on the state variable $y$, a different integrability condition of the terminal value will be required to guarantee existence and uniqueness for the solution of BSDE \eqref{BSDE:4.1}. Novel techniques and ideas will also be required to handle the unbounded growth of $g$ in $y$, which is a challenging issue. It may be further tackled in our future work.
\end{itemize}}
\end{remark}

\section{Proof of Theorem \ref{ncth1}}
\label{Proof of Theorem 1}
\setcounter{equation}{0}

This section is devoted to the proof of Theorem \ref{ncth1}. Before that, we clarify the following two basic facts. First, since $g_1$ satisfies (A1)-(A2), it follows that $f_1$ defined in \eqref{fdy1} is a convex function valued in $\R\cup\{+\infty\}$ and $\as$, for each $q\in \R^{1\times d}$, we have
\begin{align}\label{fx-nc}
f_1(\omega,t,q)\geq-\alpha_t(\omega)+\frac{1}{2\gamma}|q|^2.
\end{align}
Second, let $q$ be in $\A$, if this set is not empty. In light of (H1)-(H2) of $g_2$ and the definition of $\A$, thanks to Girsanov's theorem and Proposition \ref{L1}, BSDE \eqref{bsde-qnc} admits a unique $L^1$ solution $(Y^q_\cdot, Z^q_\cdot)$ \vspace{0.3cm}under $\Q_q$.

In the sequel, we prove Theorem \ref{ncth1} with the above observations.\vspace{0.2cm}

\noindent\textbf{Proof of (i) of Theorem \ref{ncth1}.}
Suppose first that $|\xi|+\int_0^T \alpha_t\dif t$ admits an exponential moment of $p$-order for some $p>\gamma$ (i.e., \eqref{eq:4.2}). Suppose further that $g_1(\omega,t,\cdot)$ is quadratic and convex (i.e., (A1) and (A2)), and that $g_2(\omega,t,\cdot)$ is bounded and uniformly continuous (i.e., (H1) and (H2) with $a=0$). \vspace{0.2cm}

It is not hard to verify that $g:=g_1+g_2$ also satisfies assumption (A1), hence the existence result has been given in (i) of Proposition \ref{cun}. That is to say, BSDE \eqref{BSDE:4.1}  admits a solution $(Y_\cdot, Z_\cdot)$ satisfying \eqref{kj}, and $Z_\cdot\in\M^2$. Now, we divide the following proof into three steps to show the uniqueness. \vspace{0.2cm}

\textbf{First step.} Let us start by showing $Y_\cdot\leq Y_\cdot^q$ for any $q\in\A$. Thanks to Girsanov's theorem, BSDE \eqref{BSDE:4.1}  can be equivalently written as follows:
\begin{align*}
Y_t=\xi+\int_t^T\big(q_sZ_s-g(s,Z_s)\big)\dif s+\int_t^TZ_s\cdot\dif B_s^q,\ \ t\in[0,T].
\end{align*}
It follows from \eqref{fdy1} that
\begin{align}\label{z9}
\notag \mathbbm{1}_{\{Y_t>Y^q_t\}}&\Big(q_sZ_s-g(s,Z_s)-\big(f_1(s,q_s)-g_2(s,Z_s)\big)\Big)\\
\notag =&\mathbbm{1}_{\{Y_t>Y^q_t\}}\Big(q_sZ_s-g_1(s,Z_s)-g_2(s,Z_s)-f_1(s,q_s)+g_2(s,Z_s)\Big)\\
=&\mathbbm{1}_{\{Y_t>Y^q_t\}}\Big(q_sZ_s-g_1(s,Z_s)-f_1(s,q_s)\Big)\leq0.
\end{align}
Since $g_2$ satisfies (H1)-(H2), then the generator $f_1(s,q_s)-g_2(s,\cdot)$ in BSDE \eqref{bsde-qnc} also satisfies (H1)-(H2). Thanks to (i) of Proposition \ref{bj-nc}, it suffices to prove that $(Y_\cdot-Y^q_\cdot)^+\in\s$ under $\Q_q$. Next, we will prove that $(Y_\cdot-Y^q_\cdot)^+$ is a bounded process. For each $t\in[0,T]$ and each integer $m\geq1$, we set
 $$\tau_m^t:=\inf\Big\{s\geq t: \int_t^s|Z_u|^2\dif u+\int_t^s|Z^q_u|^2\dif u\geq m\Big\}\wedge T$$
with the convention $\inf\emptyset=+\infty$. In view of $g_2$ satisfying (H1) and (H2) with $a=0$, applying It\^o-Tanaka's formula to $(Y_s-Y_s^q)^+$ and using \eqref{fdy1}, we obtain
\begin{align*}
\dif (Y_s-Y_s^q)^+&\geq \mathbbm{1}_{\{Y_s-Y_s^q>0\}}\big(g(s,Z_s)-q_sZ_s+f_1(s,q_s)-g_2(s,Z_s^q)\big)\dif s-\mathbbm{1}_{\{Y_s-Y_s^q>0\}}(Z_s-Z^q_s)\cdot\dif B_s^q\\
&\geq \mathbbm{1}_{\{Y_s-Y_s^q>0\}}\big(g_2(s,Z_s)-g_2(s,Z_s^q)\big)\dif s-\mathbbm{1}_{\{Y_s-Y_s^q>0\}}(Z_s-Z^q_s)\cdot\dif B_s^q\\
&\geq -2b ~\dif s-\mathbbm{1}_{\{Y_s-Y_s^q>0\}}(Z_s-Z^q_s)\cdot\dif B_s^q, \ \ s\in[t, \tau_m^t].
\end{align*}
Hence, we have
\begin{align*}
(Y_t-Y_t^q)^+\leq(Y_{\tau^t_m}-Y_{\tau^t_m}^q)^++2bT+\int_t^{\tau^t_m}\mathbbm{1}_{\{Y_s-Y_s^q>0\}}(Z_s-Z^q_s)\cdot\dif B_s^q.
\end{align*}
It then follows that
\begin{align}\label{a}
(Y_t-Y_t^q)^+\leq\E^{\Q_q}\big[\big(Y_{\tau^t_m}-Y_{\tau^t_m}^q\big)^+\big|\F_t\big]+2bT.
\end{align}
In view of $Y_{\tau^t_m}\leq \sup_{t\in[0,T]} |Y_t|$, the Fenchel's inequality gives
\begin{align*}
\E^{\Q_q}\Big[\sup_{t\in\T} |Y_t|\Big]&=\E\Big[M^q_T\sup_{t\in\T} |Y_t|\Big]\\
&\leq \E\Big[\exp\Big(p\sup_{t\in\T} |Y_t|\Big)\Big]+\frac{1}{p}\E\Big[M^q_T\Big(\ln M^q_T-\ln p-1\Big)\Big]\\
&=\E\Big[\exp\Big(p\sup_{t\in\T} |Y_t|\Big)\Big]+\frac{1}{p}\E\Big[M^q_T\ln M^q_T\Big]-\frac{1}{p}\Big(\ln p+1\Big).
\end{align*}
Some uncomplicated calculations give
\begin{align}\label{t}
\notag \E[M^q_T\ln M^q_T]&=\E^{\Q_q}\big[\ln {M^q_T}\big]=\E^{\Q_q}\Big[\int_0^T q_s\dif B_s-\frac{1}{2}\int_0^T|q_s|^2\dif s\Big]\\
&=\E^{\Q_q}\Big[\int_0^T q_s\dif B^q_s+\frac{1}{2}\int_0^T |q_s|^2\dif s\Big]=\frac{1}{2}\E^{\Q_q}\Big[\int_0^T |q_s|^2\dif s\Big].
\end{align}
The calculation similar to \eqref{t} will be used several times later, and this calculation process will not be repeated in detail. For ease of notations, we denote the sum of constants by $C_p$, then we have
\begin{align*}
\E^{\Q_q}\Big[\sup_{t\in\T} |Y_t|\Big]\leq C_p+\frac{1}{2p}\E^{\Q_q}\Big[\int_{0}^{T}|q_s|^2\dif s\Big]<+\infty.
\end{align*}
Furthermore, since $(Y_\cdot^q,Z_\cdot^q)$ is an $L^1$ solution under $\Q_q$, we know that $Y_\cdot^q$ belongs to class $(D)$ under $\Q_q$. Then letting $m\rightarrow\infty$ in \eqref{a} gives
\begin{align*}
(Y_t-Y_t^{q})^+\leq\E^{\Q_q}\big[(Y_{T}-Y_{T}^{q})^+\big|\F_t\big]+2bT=2bT,
\end{align*}
which naturally yields that $(Y_\cdot-Y^q_\cdot)^+\in\s$ under $\Q_q$. Hence, in view of (i) of Proposition \ref{bj-nc}, we can conclude that for \vspace{0.2cm}each $t\in\T$, $\ps$, $Y_t\leq Y_t^{q}$.

\textbf{Second step.} Set $q_s^*\in\partial_z g_1(s,Z_s)$. Then we have
$$
f_1(s,q_s^*)=q_s^*Z_s-g_1(s,Z_s), ~~s\in[0,T].
$$
Consequently, if $q^*\in\A$, then by BSDEs \eqref{BSDE:4.1}  and \eqref{bsde-qnc} along with the uniqueness of the $L^1$ solution of BSDE \eqref{bsde-qnc} under $\Q_{q^*}$, we can conclude that $Y_\cdot=Y_\cdot^{q^*}$.\vspace{0.2cm}

\textbf{Third step.} We conclude the proof by verifying $q^*\in\A$. Thanks to \eqref{fx-nc} and Young's inequality, we have
\begin{align}\label{q11}
\notag g_1(s,Z_s)&=q_s^*Z_s-f_1(s,q_s^*)\leq \frac{1}{2}\Big(\frac{1}{2\gamma}|q_s^*|^2+2\gamma|Z_s|^2\Big)+\alpha_s-\frac{1}{2\gamma}|q_s^*|^2\\
&=\gamma|Z_s|^2+\alpha_s-\frac{1}{4\gamma}|q_s^*|^2;\\
\notag \frac{1}{4\gamma}|q_s^*|^2&\leq-g_1(s,Z_s)+\gamma|Z_s|^2+\alpha_s,\ \ s\in[0,T],
\end{align}
which implies that $\ps$, $\int_0^T|q_s^*|^2\dif s<+\infty$. For each $n\geq1$, define
$$
M^*_t:=\exp\Big(\int_{0}^{t}q_s^*\dif B_s-\frac{1}{2}\int_{0}^{t}|q_s^*|^2\dif s\Big),\ \ t\in\T,
$$
$$
\tau_n:=\inf \Big\{t \in[0, T]: \int_0^t|q_s^*|^2 \dif s+\int_0^t|Z_s|^2 \dif s\geq n\Big\} \wedge T, \quad \frac{\dif \Q_n^*}{\dif \mathbb{P}}:=M^*_{\tau_n}.
$$
We set $B_t^{q^*}:=B_t-\int_0^t (q^*_s)^{\top}\dif s, ~t\in\T$, then $(B_t^{q^*})_{t\in[0, \tau_n]}$ is a standard $d$-dimensional Brownian motion under the probability measure $\Q_n^*$ for each $n\geq1$. Now, let us further show that $(M^*_{\tau_n})_n$ is uniformly integrable.\vspace{0.2cm}

\begin{lemma}\label{lem1}
$(M^*_{\tau_n})_n$ is uniformly integrable.\vspace{0.2cm}
\end{lemma}

\noindent\textbf{Proof.} In view of \eqref{kj} and the Fenchel's inequality we have
\begin{align*}
\E^{\Q_n^{*}}\Big[\sup_{t\in\T} \Big(|Y_t|+\int_0^t\alpha_s\dif s\Big)\Big]&=\E\Big[M^*_{\tau_n} \sup_{t\in\T} \Big(|Y_t|+\int_0^t\alpha_s\dif s\Big)\Big]\\
&\leq C_p+\frac{1}{2p} \E^{\Q_n^*}\Big[\int_0^{\tau_n}|q_s^*|^2 \dif s\Big].
\end{align*}
Since $g_1(s,Z_s)=q_s^*Z_s-f_1(s,q_s^*)$, by BSDE \eqref{BSDE:4.1}  we have
$$
Y_0=Y_{\tau_n}+\int_0^{\tau_n}\big(f_1(s,q_s^*)-g_2(s,Z_s)\big)\dif s+\int_0^{\tau_n} Z_s \cdot \dif B_s^{q^*}.
$$
By virtue of assumption (H2) with $a=0$, we obtain that $\as$, for each $z\in\R^{d}$, $|g_2(\omega,t,z)|\leq b$. Hence, in view of \eqref{fx-nc} and the fact that $Y_{\tau_n}\geq-|Y_{\tau_n}|$, we can deduce that
\begin{align*}
Y_0&=Y_{\tau_n}+\int_0^{\tau_n}\big(f_1(s,q_s^*)-g_2(s,Z_s)\big)\dif s+\int_0^{\tau_n} Z_s \cdot \dif B_s^{q^*}\\
&\geq-|Y_{\tau_n}|-\int_0^{\tau_n}\alpha_s \dif s+\frac{1}{2\gamma}\int_0^{\tau_n}|q_s^*|^2\dif s-\int_0^{\tau_n}b ~\dif s+\int_0^{\tau_n}Z_s \cdot \dif B_s^{q^*}\\
&\geq-\sup_{t\in\T}\Big(|Y_t|+\int_0^t\alpha_s\dif s\Big)+\frac{1}{2\gamma}\int_0^{\tau_n}|q_s^*|^2\dif s-bT+\int_0^{\tau_n}Z_s \cdot \dif B_s^{q^*}.
\end{align*}
It then follows that
\begin{align}\label{p}
\notag Y_0&\geq-\E^{\Q_n^*}\Big[\sup_{t\in\T}\Big(|Y_t|+\int_0^t\alpha_s\dif s\Big)\Big]+\frac{1}{2\gamma}\E^{\Q_n^*}\Big[\int_0^{\tau_n}|q_s^*|^2\dif s\Big]-bT\\
&\geq-\big(C_p+bT\big)+\frac{1}{2}\Big(\frac{1}{\gamma}-\frac{1}{p}\Big)\E^{\Q_n^*}\Big[\int_0^{\tau_n}|q_s^*|^2\dif s\Big].
\end{align}
Since $p>\gamma$, then $\frac{1}{2}\Big(\frac{1}{\gamma}-\frac{1}{p}\Big)>0.$ Hence, in view of \eqref{t} and \eqref{p}, we obtain
\begin{align}\label{g}
2\E[M^*_{\tau_n}\ln M^*_{\tau_n}]=2\E^{\Q_n^{*}}[\ln M^*_{\tau_n}]=\E^{\Q_n^{*}}\Big[\int_0^{\tau_n}|q_s^*|^2\dif s\Big]\leq C_{p,\gamma,T,b},
\end{align}
where $C_{p,\gamma,T,b}$ is a positive constant independent of $n$. Hence, we obtain $\sup_n\E[M^*_{\tau_n}\ln M^*_{\tau_n}]<+\infty$. Then $(M^*_{\tau_n})_n$ is uniformly integrable by Lemma \ref{lem-poussin} (the de La Vall$\acute{\rm{e}}$e Poussin lemma), and $\E[M^*_T]=1$, which implies that $\mathcal{E}(q^*):=(M^*_t)_{t\in\T}$ is a uniformly integrable martingale and defines a probability measure $\Q^*$ by $\frac{\dif \Q^*}{\dif \mathbb{P}}:=M^*_{T}$. The proof of Lemma \ref{lem1} is then \vspace{0.3cm}complete. \hfill\framebox

Moreover, applying Fatou's lemma in \eqref{g}, we have
$$
2\E[M^*_T\ln M^*_T]=\E^{\Q^*}\Big[\int_0^T|q_s^*|^2\dif s\Big]\leq\liminf_n\E^{\Q_n^* }\Big[\int_0^{\tau_n}|q_s^*|^2\dif s\Big]<+\infty.
$$
Thus, it remains to show that
\begin{align}\label{e}
\E^{\Q^*}\Big[\int_0^T|f_1(s,q_s^*)|\dif s\Big]<+\infty.
\end{align}
In view of \eqref{fx-nc}, we know that
$$
f_1^-(s,q_s^*)\leq\alpha_s,
$$
so, we have
\begin{align}\label{fb}
\E^{\Q^*}\Big[\int_0^T f_1^-(s,q_s^*)\dif s\Big]\leq\E^{\Q^*}\Big[\int_0^T \alpha_s\dif s\Big]<+\infty.
\end{align}
Moreover, thanks to BSDE \eqref{BSDE:4.1} we have for each $n\geq 1$,
\begin{align*}
Y_0&=Y_{\tau_n}+\int_0^{\tau_n} \big(f_1(s,q_s^*)-g_2(s,Z_s)\big)\dif s+\int_0^{\tau_n} Z_s\cdot\dif B_s^{q^*}\\
&=\E^{\Q^*}\Big[Y_{\tau_n}+\int_0^{\tau_n}\big(f_1(s,q_s^*)-g_2(s,Z_s)\big)\dif s\Big].
\end{align*}
It then follows that
\begin{align*}
Y_0&\geq\E^{\Q^*}\Big[-|Y_{\tau_n}|+\int_0^{\tau_n} f_1^+(s,q_s^*)\dif s-\int_0^{\tau_n} f^-(s,q_s^*)\dif s-\int_0^{\tau_n} g_2(s,Z_s)\dif s\Big]\\
&\geq\E^{\Q^*}\Big[-|Y_{\tau_n}|+\int_0^{\tau_n} f_1^+(s,q_s^*)\dif s-\int_0^{\tau_n} f_1^-(s,q_s^*)\dif s\Big]-bT.
\end{align*}
Thus, for each $n\geq 1$, we have
\begin{align}\label{zb}
\E^{\Q^*}\Big[\int_0^{\tau_n} f_1^+(s,q_s^*)\dif s\Big]\leq Y_0+\E^{\Q^*}\Big[\sup\limits_{t\in\T}|Y_t|\Big]+\E^{\Q^*}\Big[\int_0^T f_1^-(s,q_s^*)\dif s\Big]+bT<+\infty.
\end{align}
Combining \eqref{fb} and \eqref{zb} along with Fatou's lemma yields \eqref{e}. Thus, we have shown that $q^*\in\A$ and $Y_\cdot=Y^{q^*}_\cdot=\essinf_{q\in\A}Y^{q}_\cdot$, which yields the uniqueness result. The proof of (i) of Theorem \ref{ncth1} is complete. \vspace{0.3cm}\hfill\framebox

\noindent\textbf{Proof of (ii) of Theorem \ref{ncth1}.}
Suppose that $|\xi|+\int_0^T \alpha_t\dif t$ admits an exponential moment of $p$-order for some $p>\gamma$ (i.e., \eqref{eq:4.2}). Suppose further that $g_1(\omega,t,\cdot)$ is convex with a strict quadratic growth (i.e., (A1') and (A2)), and that $g_2(\omega,t,\cdot)$ is uniformly continuous with a sub-linear growth (i.e., (H1)-(H2)). \vspace{0.2cm}

It is not difficult to verify that the generator $g:=g_1+g_2$ of BSDE \eqref{BSDE:4.1} also satisfies assumption (A1') with constant $\bar\gamma$ and $\gamma$ being respectively replaced with $\bar\gamma-\eps$ and $\gamma+\eps$ for any constant $\eps\in (0,\bar\gamma)$. Note that (A1') implies (A1) and that $p>\gamma+\eps$ when $\eps$ is small enough. It follows from (i) of Proposition \ref{cun} along with Proposition \ref{Proposition 2 in [FanHuTang 2020]} that BSDE \eqref{BSDE:4.1} admits a solution $(Y_\cdot,Z_\cdot)$ satisfying \eqref{kj}, \eqref{z1} and \eqref{z2}. For the uniqueness, we will fit the proof of (i) of Theorem \ref{ncth1} to this situation. The following proof will be divided into three steps.\vspace{0.2cm}

\textbf{First step.} Let us start with showing $Y_\cdot\leq Y_\cdot^q$ for any $q\in\A$. Thanks to Girsanov's theorem, BSDE \eqref{BSDE:4.1}  can be equivalently written as follows:
\begin{align}\label{4.10xing}
Y_t=\xi+\int_t^T\big(q_sZ_s-g(s,Z_s)\big)\dif s+\int_t^TZ_s\cdot\dif B_s^q,\ \ t\in[0,T].
\end{align}
In view of \eqref{kj} and Fenchel's inequality, we have
\begin{align*}
\E^{\Q_q}\Big[\sup_{t\in\T} |Y_t|\Big]&=\E\Big[M^q_T \sup_{t\in\T} |Y_t|\Big]\leq C_p+\frac{1}{2p} \E^{\Q_q}\Big[\int_0^T|q_s|^2 \dif s\Big]<+\infty,
\end{align*}
and in the same manner, by virtue of \eqref{z1} we have
\begin{align*}
\E^{\Q_q}\Big[\int_0^T |Z_t|^2\dif s\Big]\leq C_\eta+\frac{1}{2\eta} \E^{\Q_q}\Big[\int_0^T|q_s|^2 \dif s\Big]<+\infty.
\end{align*}
Consequently, $(Y_\cdot,Z_\cdot)$ is an $L^1$ solution of BSDE \eqref{4.10xing} under $\Q_q$. Note that $g_2$ satisfies assumptions (H1)-(H2), $(Y^q_\cdot, Z^q_\cdot)$ is an $L^1$ solution of BSDE \eqref{bsde-qnc} under $\Q_q$ and that \eqref{z9} still holds. We can apply (ii) of Proposition \ref{bj-nc} to obtain that for each $t\in\T$, $\ps$, $Y_t\leq Y^q_t$.\vspace{0.3cm}

\textbf{Second step.} Set $q_s^*\in\partial_z g_1(s,Z_s)$. Then we have
\begin{align}\label{4.11xing-2}
f_1(s,q_s^*)=q_s^*Z_s-g_1(s,Z_s), ~~s\in[0,T].
\end{align}
Consequently, if $q^*\in\A$, then by BSDEs \eqref{BSDE:4.1}  and \eqref{bsde-qnc} along with the uniqueness of the $L^1$ solution of BSDE \eqref{bsde-qnc} under $\Q_{q^*}$, we can conclude that $Y_\cdot=Y_\cdot^{q^*}$.\vspace{0.3cm}

\textbf{Third step.} We conclude the proof by verifying $q^*\in\A$. By \eqref{q11} we have
$$\ps,\ \ \int_0^T|q_s^*|^2\dif s<+\infty.$$
For each $n\geq1$, let us define
$$
M^*_t:=\exp\Big(\int_{0}^{t}q_s^*\dif B_s-\frac{1}{2}\int_{0}^{t}|q_s^*|^2\dif s\Big),\ \ t\in\T,
$$
$$
\tau_n:=\inf \Big\{t \in[0, T]: \int_0^t|q_s^*|^2 \dif s+\int_0^t|Z_s|^2 \dif s\geq n\Big\} \wedge T, \quad \frac{\dif \Q_n^*}{\dif \mathbb{P}}:=M^*_{\tau_n}.
$$
We set
$$B_t^{q^*}:=B_t-\int_0^t (q^*_s)^{\top}\dif s,\ \ t\in\T,$$
then $(B_t^{q^*})_{t\in[0,\tau_n]}$ is a standard $d$-dimensional Brownian motion under the probability measure $\Q_n^*$ for each $n\geq1$. Let us further show that $\E^{\Q^*}\Big[\int_{0}^{T}|q^*_s|^2\dif s\Big]<+\infty$, and $\big(M^*_t\big)_{t\in[0,T]}$ is a uniformly integrable martingale. To do this, we will show that Lemma \ref{lem1} still holds in such case.\vspace{0.2cm}

As discussed in Lemma \ref{lem1}, we have
\begin{align}\label{4.11xing}
\E^{\Q_n^{*}}\Big[\sup_{t\in\T}\Big(|Y_t|+\int_0^t\alpha_s\dif s\Big)\Big]\leq C_p+\frac{1}{2p} \E^{\Q_n^*}\Big[\int_0^{\tau_n}|q_s^*|^2 \dif s\Big].
\end{align}
Moreover, due to $\iota\in(0,1)$, it is easy to verify that
$$
\int_0^{\tau_n}|Z_s|^\iota\dif s\leq T+\int_0^T|Z_s|\dif s.
$$
Taking \eqref{z2} and \eqref{t} into consideration, for each $\lambda>0$, we have
\begin{align}\label{zz}
\notag \E^{\Q_n^{*}}\Big[\int_0^{\tau_n}|Z_s|^\iota\dif s\Big]
&=\E\Big[M^*_{\tau_n}\int_0^{\tau_n}|Z_s|^\iota\dif s\Big] \\
\notag &\leq \E\Big[\exp\Big(\lambda\int_0^{\tau_n}|Z_s|^\iota\dif s\Big)\Big]+\frac{1}{\lambda}\E\big[M^*_{\tau_n}(\ln M^*_{\tau_n }-\ln \lambda-1)\big] \\
\notag &\leq \E\Big[\exp\Big(T\lambda+\lambda\int_0^T|Z_s|\dif s\Big)\Big]+\frac{1}{\lambda}\E\big[M^*_{\tau_n}(\ln M^*_{\tau_n }-\ln \lambda-1)\big]\\
&\leq C_{T, \lambda}+\frac{1}{2\lambda} \E^{\Q_n^*}\Big[\int_0^{\tau_n}|q_s^*|^2 \dif s\Big],
\end{align}
where $C_{T, \lambda}$ is a positive constant depending only on $(T, \lambda)$. Thus, in view of assumption (H2), \eqref{BSDE:4.1} , \eqref{4.11xing-2}, \eqref{fx-nc} and the fact that $Y_{\tau_n}\geq-|Y_{\tau_n}|$, we have
\begin{align*}
Y_0&=Y_{\tau_n}+\int_0^{\tau_n}\big(f_1(s,q_s^*)-g_2(s,Z_s)\big)\dif s+\int_0^{\tau_n} Z_s \cdot \dif B_s^{q^*}\\
&\geq-|Y_{\tau_n}|-\int_0^{\tau_n}\alpha_s \dif s+\frac{1}{2\gamma}\int_0^{\tau_n}|q_s^*|^2\dif s-\int_0^{\tau_n}(a|Z_s|^\iota+b)\dif s+\int_0^{\tau_n}Z_s \cdot \dif B_s^{q^*}\\
&\geq-\sup_{t\in\T}\Big(|Y_t|+\int_0^t\alpha_s\dif s\Big)+\frac{1}{2\gamma}\int_0^{\tau_n}|q_s^*|^2\dif s-a\int_0^{\tau_n}|Z_s|^\iota\dif s-bT+\int_0^{\tau_n}Z_s \cdot \dif B_s^{q^*}.
\end{align*}
It then follows that, in view of \eqref{4.11xing} and \eqref{zz},
\begin{align*}
Y_0&\geq-\E^{\Q_n^*}\Big[\sup_{t\in\T}\Big(|Y_t|+\int_0^t\alpha_s\dif s\Big)\Big]+\frac{1}{2\gamma}\E^{\Q_n^*}\Big[\int_0^{\tau_n}|q_s^*|^2\dif s\Big]-a\E^{\Q_n^*}\Big[\int_0^{\tau_n}|Z_s|^\iota\dif s\Big]-bT \\
&\geq-C_p-aC_{T, \lambda}-bT +\frac{1}{2}\Big(\frac{1}{\gamma}-\frac{1}{p}-\frac{a}{\lambda}\Big)\E^{\Q_n^*}\Big[\int_0^{\tau_n}|q_s^*|^2\dif s\Big].
\end{align*}
Since $p>\gamma$, setting $\lambda>\frac{ap\gamma}{p-\gamma}$, then we have $\frac{1}{2}\Big(\frac{1}{\gamma}-\frac{1}{p}-\frac{a}{\lambda}\Big)>0.$ Hence,
\begin{align}\label{49}
2\E[M^*_{\tau_n}\ln M^*_{\tau_n}]=2\E^{\Q_n^{*}}[\ln M^*_{\tau_n}]=\E^{\Q_n^{*}}\Big[\int_0^{\tau_n}|q_s^*|^2\dif s\Big]<C_{p,\gamma,T,a,b,\lambda},
\end{align}
where $C_{p,\gamma,T,a,b,\lambda}>0$ is a constant independent of $n$. Thus, $\sup_n\E[M^*_{\tau_n}\ln M^*_{\tau_n}]<+\infty$. Thanks to Lemma \ref{lem-poussin} (the de La Vall$\acute{\rm{e}}$e Poussin lemma), we obtain that $(M^*_{\tau_n})_n$ is uniformly integrable, and $\E[M^*_T]=1$, which implies that $\mathcal{E}(q^*):=(M^*_t)_{t\in\T}$ is a uniformly integrable martingale and defines a probability measure $\Q^*$ by $\frac{\dif \Q^*}{\dif \mathbb{P}}:=M^*_{T}$. That is to say, Lemma \ref{lem1} holds still \vspace{0.3cm}in this case.

Moreover, applying Fatou's lemma in \eqref{49}, we have
\begin{align}\label{4.12xing}
2\E[M^*_T\ln M^*_T]=\E^{\Q^*}\Big[\int_0^T|q_s^*|^2\dif s\Big]\leq\liminf_n\E^{\Q_n^* }\Big[\int_0^{\tau_n}|q_s^*|^2\dif s\Big]<+\infty.
\end{align}
Hence, it remains to show that $\E^{\Q^*}\big[\int_0^T|f_1(s,q_s^*)|\dif s\big]<+\infty$. In view of \eqref{fx-nc}, we deduce that
$$
f_1^-(s,q_s^*)\leq\alpha_s,
$$
which means that \eqref{fb} still holds. Moreover, thanks to BSDE \eqref{BSDE:4.1} we have for each $n\geq 1$,
\begin{align*}
Y_0&=Y_{\tau_n}+\int_0^{\tau_n} \big(f_1(s,q_s^*)-g_2(s,Z_s)\big)\dif s+\int_0^{\tau_n} Z_s\cdot\dif B_s^{q^*}\\ &=\E^{\Q^*}\Big[Y_{\tau_n}+\int_0^{\tau_n}\big(f_1(s,q_s^*)-g_2(s,Z_s)\big)\dif s\Big].
\end{align*}
By virtue of (H2), it then follows that
\begin{align*}
Y_0\geq\E^{\Q^*}\Big[-|Y_{\tau_n}|+\int_0^{\tau_n} f_1^+(s,q_s^*)\dif s-\int_0^{\tau_n} f_1^-(s,q_s^*)\dif s-a\int_0^{\tau_n} |Z_s|^\iota\dif s-bT\Big].
\end{align*}
Thus, by Fenchel's inequality, \eqref{kj}, \eqref{z2}, \eqref{fb}, \eqref{zz} and \eqref{4.12xing} we have for each $n\geq 1$,
\begin{align*}
\E^{\Q^*}\Big[\int_0^{\tau_n} f_1^+(s,q_s^*)\dif s\Big]\leq &Y_0+\E^{\Q^*}\Big[\sup\limits_{t\in\T}|Y_t|\Big]+\E^{\Q^*}\Big[\int_0^T f_1^-(s,q_s^*)\dif s\Big]\\
&+a\E^{\Q^*}\Big[\int_0^T |Z_s|\dif s\Big]+bT<+\infty.
\end{align*}
Combining \eqref{fb} and the last inequality along with Fatou's lemma yields
$$
\E^{\Q^*}\Big[\int_0^{\tau_n}|f_1(s,q_s^*)|\dif s\Big]=\E^{\Q^*}\Big[\int_0^T\big(f_1^+(s,q_s^*)+f_1^-(s,q_s^*)\big)\dif s\Big]<+\infty.
$$
Thus, we have shown that $q^*\in\A$ is optimal, i.e., $Y_\cdot=Y^{q^*}_\cdot=\essinf_{q\in\A}Y^{q}_\cdot$, which naturally yields the uniqueness result. The proof of (ii) of Theorem \ref{ncth1} is then \vspace{0.3cm} complete.\hfill\framebox

\noindent\textbf{Proof of (iii) of Theorem \ref{ncth1}.}
Suppose that $|\xi|+\int_0^T \alpha_t\dif t$ admits an exponential moment of $p$-order for some $p>\gamma$ (i.e., \eqref{eq:4.2}), and that $g_1(\omega,t,\cdot)$ is quadratic and strongly convex (i.e., (A1) and (A2')), and $g_2(\omega,t,\cdot)$ is uniformly continuous (i.e., (H1)). Moreover, we also suppose that there exists a process $u_t(\omega)\in\partial_z g_1(\omega,t,0)$ such that $\as$, $|u_t(\omega)|^2\leq\alpha_t(\omega)$. \vspace{0.2cm}

According to Proposition \ref{observation 1}, we know that $g_1$ also satisfies assumption (A1'). Since $g_2$ satisfies (H1), it follows that the generator $g:=g_1+g_2$ of BSDE \eqref{BSDE:4.1} also satisfies (A1') with constant $\bar\gamma$ and $\gamma$ being respectively replaced with $\bar\gamma-\eps$ and $\gamma+\eps$ for any constant $\eps\in (0,\bar\gamma)$. Note that (A1') implies (A1) and that $p>\gamma+\eps$ when $\eps$ is small enough. It follows from (i) of Proposition \ref{cun} along with Proposition \ref{Proposition 2 in [FanHuTang 2020]} that BSDE \eqref{BSDE:4.1} admits a solution $(Y_\cdot,Z_\cdot)$ satisfying \eqref{kj}, \eqref{z1} and \eqref{z2}. \vspace{0.2cm}

For the uniqueness, we first present the following proposition on the uniform integrability of the solution. It can be proved identically as in the third step of the proof of (ii) of Theorem \ref{ncth1}. Here we omit the details.

\begin{proposition}\label{pro-Q-nc-111}
Assume that $g:=g_1+g_2$ is a generator such that $g_1$ satisfies (A1') and (A2), and $g_2$ satisfies (H1). Let $(Y_\cdot, Z_\cdot)$ be a solution of BSDE \eqref{BSDE:4.1}  such that  \eqref{kj}, \eqref{z1} and \eqref{z2} hold for some $\eta>0$ depending only on $(\gamma, \bar {\gamma}, T, p)$ and each $\lambda>0$. Then, for all $(\F_t)$-progressively measurable process $(q^*_s)_{s\in[0,T]}$ valued in $\R^{1\times d}$ and such that $q_s^*\in\partial_z g_1(s,Z_s)$ for all $s\in[0,T]$, $\mathcal{E}(q^*)$ is a uniformly integrable martingale and defines a probability measure $\Q^*\sim\mathbb{P}$. Moreover, we have
\begin{align*}
\E^{\Q^*}\Big[\int_0^T|q^*_s|^2\dif s\Big]<+\infty.
\end{align*}
\end{proposition}

Next, let $(Y_\cdot, Z_\cdot)$ and $(Y'_\cdot,Z'_\cdot)$ be two solutions of BSDE \eqref{BSDE:4.1}  such that both $(Y_\cdot, Z_\cdot)$ and $(Y'_\cdot,Z'_\cdot)$ satisfy \eqref{kj}, \eqref{z1} and \eqref{z2} for $\eta>0$ depending only on $(\gamma, \bar {\gamma}, T, p)$ and each $\lambda>0$. In order to verify the uniqueness, by a symmetry argument it is sufficient to show that $\ps, ~Y_t\geq Y'_t$ for each $t\in[0,T)$. For $t\in[0,T)$, let us denote $D:=\{Y_t<Y'_t\}$, and set the stopping time $\tau:=\inf\{s\geq t | Y_s\geq Y'_s\}$. Then for each $s\in [t,\tau]$, we have $\ps, ~Y_s\mathbbm{1}_D \leq Y'_s\mathbbm{1}_D $ and $Y_\tau\mathbbm{1}_D = Y'_\tau\mathbbm{1}_D$ since $\ps$, $t\rightarrow Y_t$ is continuous. It implies that $((Y_s-Y'_s)\mathbbm{1}_D )_{s\in[t,\tau]}$ is a non-positive process. In the sequel, we will first prove that it is a bounded process.\vspace{0.2cm}

Let us consider an $(\F_t)$-progressively measurable process $(q_s^*)_{s\in[0,T]}$ valued in $\R^{1\times d}$ and such that $q_s^*\in\partial_z g_1(s,Z_s)$ for all $s\in[0,T]$. According to Proposition \ref{pro-Q-nc-111} and Girsanov's theorem, we know that $\mathcal{E}(q^*)$ defines a probability measure $\Q^*$ and $B^{q^*}_t:=B_t-\int_0^t (q^*_s)^{\top}\dif s, ~t\in[0,T]$ is a standard $d$-dimensional Brownian motion under the probability measure $\Q^*$. Hence, in view of (H1) of $g_2$, we have
\begin{align}\label{dd-nc-222}
\notag \dif \big(Y_s-Y'_s\big)&=\Big(g(s,Z_s)-g(s,Z'_s)-q_s^*(Z_s-Z'_s)\Big)\dif s-(Z_s-Z'_s)\cdot\dif B_s^{q^*}\\
\notag &\leq\Big(g_1(s,Z_s)-g_1(s,Z'_s)-q_s^*(Z_s-Z'_s)+\phi(|Z_s-Z'_s|)\Big)\dif s-(Z_s-Z'_s)\cdot\dif B_s^{q^*}\\
\notag &\leq\Big(g_1(s,Z_s)-g_1(s,Z'_s)-q_s^*(Z_s-Z'_s)+A|Z_s-Z'_s|+A\Big)\dif s\\
&\ \ \ \ -(Z_s-Z'_s)\cdot\dif B_s^{q^*}, \ \ \ s\in\T.
\end{align}
Let $q_s^A:=\frac{A(Z_s-Z'_s)}{|Z_s-Z'_s|}\mathbbm{1}_{|Z_s-Z'_s|\neq0}$, and define the probability measure $\Q^A$ equivalent to $\Q^*$ by
\begin{align*}
\frac{\dif \Q^A}{\dif \Q^*}:=\exp\left(\int_0^T q_s^A\cdot\dif B^{q^*}_s-\frac{1}{2}\int_0^T |q_s^A|^2\dif s\right).
\end{align*}
We set $B^A_t:=B^{q^*}_t-\int_0^t q^A_s\dif s, ~t\in[0,T]$, then Girsanov's theorem gives that $(B^A_t)_{t\in\T}$ is a standard $d$-dimensional Brownian motion under $\Q^A$. It then follows from \eqref{dd-nc-222} that
\begin{align*}
\dif \big(Y_s-Y'_s\big)&\leq\Big(g_1(s,Z_s)-g_1(s,Z'_s)-q_s^*(Z_s-Z'_s)+A\Big)\dif s-(Z_s-Z'_s)\cdot\dif B_s^A, \ \ s\in\T.
\end{align*}
Recall that $c$ and $\varepsilon$ are defined in assumption (A2') of $g_1$, and set $\bar{c}:=c+A$. Applying It\^o's formula to $e^{\varepsilon(Y_s-Y'_s-\bar{c}s)\mathbbm{1}_D}$ and using the last inequality we can deduce that
\begin{align*}
\dif e^{\varepsilon(Y_s-Y'_s-\bar{c}s)\mathbbm{1}_D}&\leq\varepsilon\mathbbm{1}_D e^{\varepsilon(Y_s-Y'_s-\bar{c}s)\mathbbm{1}_D}\Big[-\bar{c}+\frac{1}{2}\varepsilon|Z_s-Z'_s|^2\\
&\ \ \ +\big(g_1(s,Z_s)-g_1(s,Z'_s)-q_s^*(Z_s-Z'_s)+A\big)\Big]\dif s\\
&\ \ \ -\varepsilon\mathbbm{1}_D e^{\varepsilon(Y_s-Y'_s-\bar{c}s)\mathbbm{1}_D}(Z_s-Z'_s)\cdot\dif B^A_s, \ \ s\in[t,\tau].
\end{align*}
Moreover, by virtue of assumption (A2') of $g_1$, we have
$$
-\bar{c}+\frac{1}{2}\varepsilon|Z_s-Z'_s|^2+\big(g_1(s,Z_s)-g_1(s,Z'_s)-q_s^*(Z_s-Z'_s)+A\big)\leq0,
$$
which means that
\begin{align}\label{10-nc-222}
\dif e^{\varepsilon(Y_s-Y'_s-\bar{c}s)\mathbbm{1}_D}\leq-\varepsilon\mathbbm{1}_D e^{\varepsilon(Y_s-Y'_s-\bar{c}s)\mathbbm{1}_D}(Z_s-Z'_s)\cdot\dif B^A_s,\ \ s\in[t,\tau].
\end{align}
Define
$$
\sigma^t_m:=\inf \Big\{s\geq t: \int_t^s|Z_u|^2 \dif u+\int_t^s|Z'_u|^2 \dif u\geq m\Big\} \wedge \tau.
$$
It follows from \eqref{10-nc-222} that
\begin{align*}
e^{\varepsilon(Y_s-Y'_s-\bar{c}s)\mathbbm{1}_D}\geq e^{\varepsilon(Y_{\sigma^t_m}-Y'_{\sigma^t_m}-\bar{c}\sigma^t_m)\mathbbm{1}_D}+\int_s^{\sigma^t_m}\varepsilon\mathbbm{1}_D e^{\varepsilon(Y_u-Y'_u-\bar{c}u)\mathbbm{1}_D}(Z_u-Z'_u)\cdot\dif B^A_u,\ \ s\in[t,\tau].
\end{align*}
Hence,
\begin{align*}
e^{\varepsilon(Y_s-Y'_s-\bar{c}s)\mathbbm{1}_D}
\geq\E^{\Q^A}\Big[e^{\varepsilon\big(Y_{\sigma^t_m}-Y'_{\sigma^t_m}-\bar{c}\sigma^t_m\big)\mathbbm{1}_D}\Big|\F_s\Big], \ \ s\in[t,\tau].
\end{align*}
Note that $((Y_s-Y'_s)\mathbbm{1}_D)_{s\in[t,\tau]}$ is a non-positive process. Letting $m\rightarrow\infty$ and applying Lebesgue's dominated convergence theorem to the last inequality yields
\begin{align*}
e^{\varepsilon(Y_s-Y'_s-\bar{c}s)\mathbbm{1}_D}
\geq\E^{\Q^A}\Big[e^{\varepsilon(Y_\tau-Y'_\tau-\bar{c}\tau)\mathbbm{1}_D}\big|\F_s\Big]\geq e^{-\varepsilon\bar{c}T},\ \ s\in[t,\tau].
\end{align*}
Thus, we have
$$
(Y_s-Y'_s)\mathbbm{1}_D\geq\bar{c}s-\bar{c}T\geq-(c+A)T,\ \ s\in[t,\tau],\vspace{0.2cm}
$$
which implies that $((Y_s-Y'_s)\mathbbm{1}_D )_{s\in[t,\tau]}$ is a bounded \vspace{0.3cm}non-positive process.

Furthermore, according to Lemma 4 of \cite{Fan-Jiang2011}, in light of assumptions of the function $\phi(\cdot)$ in (H1), we can conclude that for all $x\in\R^+$ and each $n\geq1$, it holds that
\begin{align}\label{eq:5.19}
\phi(x)\leq(n+2A)x+\phi(\frac{2A}{n+2A}).
\end{align}
In view of $q_s^*\in\partial_z g_1(s,Z_s)$, we have
\begin{align*}
g_1(s,Z_s)-g_1(s,Z'_s)-q_s^*(Z_s-Z'_s)\leq0,
\end{align*}
which combining \eqref{dd-nc-222} and \eqref{eq:5.19} gives
\begin{align}\label{ddn-nc-333}
\notag \dif \big(Y_s-Y'_s\big)&\leq\Big(g_1(s,Z_s)-g_1(s,Z'_s)-q_s^*(Z_s-Z'_s)+\phi(|Z_s-Z'_s|)\Big)\dif s-(Z_s-Z'_s)\cdot\dif B_s^{q^*}\\
&\leq\Big((n+2A)|Z_s-Z'_s|+\phi\Big(\frac{2A}{n+2A}\Big)\Big)\dif s-(Z_s-Z'_s)\cdot\dif B_s^{q^*}, ~~~~~s\in\T.
\end{align}
For each $n\geq1$, set
$$q^n_s:=\frac{(n+2A)(Z_s-Z'_s)}{|Z_s-Z'_s|}\mathbbm{1}_{|Z_s-Z'_s|\neq0},\ \ s\in\T,\vspace{0.2cm}$$
and define the probability measure $\Q^n$ equivalent to $\Q^*$ by
\begin{align*}
\frac{\dif \Q^n}{\dif \Q^*}:=\exp\left(\int_0^T q_s^n\cdot\dif B^{q^*}_s-\frac{1}{2}\int_0^T |q_s^n|^2\dif s\right).
\end{align*}
We set $B^n_s:=B^{q^*}_s-\int_0^s q^n_u\dif u, ~s\in[0,T]$. By Girsanov's theorem we know that $(B^n_s)_{s\in\T}$ is a standard $d$-dimensional Brownian motion under $\Q^n$ for each $n\geq1$. It then follows from \eqref{ddn-nc-333} that
\begin{align*}
\dif \big(Y_s-Y'_s\big)&\leq\phi\Big(\frac{2A}{n+2A}\Big)\dif s-(Z_s-Z'_s)\cdot\dif B_s^n,\ \ s\in\T.
\end{align*}
Define
$$
\tau^t_m:=\inf \Big\{s\geq t: \int_t^s|Z_u-Z'_u|^2 \dif u\geq m\Big\} \wedge \tau.\vspace{0.1cm}
$$
Taking integral from $t$ to $\tau_m^t$ and then the conditional expectation under $\Q^n$ in the last inequality, we obtain that for each $m, n\geq1$,
\begin{align}\label{222-nc}
\big(Y_t-Y'_t\big)\mathbbm{1}_D\geq\E^{\Q^n}\big[\big(Y_{\tau^t_m}-Y'_{\tau^t_m}\big)\mathbbm{1}_D|\F_t\big]-\phi\Big(\frac{2A}{n+2A}\Big)T\mathbbm{1}_D.
\end{align}
Since $((Y_s-Y'_s)\mathbbm{1}_D )_{s\in[t,\tau]}$ is a bounded non-positive process and $\ps$,  $Y_\tau\mathbbm{1}_D = Y'_\tau\mathbbm{1}_D$, letting $m\rightarrow\infty$ in \eqref{222-nc} and applying Lebesgue's dominated convergence theorem yields that for each $n\geq1$,
\begin{align}\label{eq:5.22}
\big(Y_t-Y'_t\big)\mathbbm{1}_D\geq-\phi\Big(\frac{2A}{n+2A}\Big)T\mathbbm{1}_D.
\end{align}
Letting $n\rightarrow\infty$ gives that
$$\phi\big(\frac{2A}{n+2A}\big) \rightarrow 0.\vspace{0.2cm}$$
Then, from inequality \eqref{eq:5.22} we deduce that $\ps$, $\big(Y_t-Y'_t\big)\mathbbm{1}_D\geq 0$. Finally, in light of the definition of $D$, we know that $\mathbb{P}(D)=0$ and then $\ps$, $Y_t\geq Y'_t$. The proof of (iii) of Theorem \ref{ncth1} is then \vspace{0.3cm}complete. \hfill\framebox

\section{Proof of Theorem \ref{ncth2}}
\label{Proof of Theorem 2}
\setcounter{equation}{0}

This section is devoted to the proof of Theorem \ref{ncth2}. Suppose that $|\xi|+\int_0^T \alpha_t\dif t$ admits an exponential moment of $\gamma$-order (i.e., \eqref{eq:4.2} holds for $p=\gamma$). Suppose further that $g_1(\omega,t,\cdot)$ is quadratic and strongly convex (i.e., (A1) and (A2')), and that $g_2(\omega,t,\cdot)$ is bounded and uniformly continuous (i.e., (H1) and (H2) with $a=0$).\vspace{0.2cm}

We first establish the following uniform integrability result.

\begin{proposition}\label{Q}
Assume that $\xi$ is a terminal value, $g=g_1+g_2$ is a generator such that $g_1$ satisfies (A1) and (A2'), $g_2$ satisfies (H2) with $a=0$, and
$$\E\Big[\exp(\gamma(|\xi|+\int_0^T\alpha_t\dif t))\Big]<+\infty.$$
Let $(Y_\cdot, Z_\cdot)$ be a solution of BSDE \eqref{BSDE:4.1}  such that $(e^{\gamma (|Y_t|+\int_0^t\alpha_s\dif s)})_{t\in[0,T]}$ belongs to class $(D)$. Then, for all $(\F_t)$-progressively measurable process $(q^*_s)_{s\in[0,T]}$ valued in $\R^{1\times d}$ and such that $q_s^*\in\partial_z g_1(s,Z_s)$ for all $s\in[0,T]$, $\mathcal{E}(q^*)$ is a uniformly integrable martingale and defines a probability measure $\Q^*$ equivalent to $\mathbb{P}$.\vspace{0.2cm}
\end{proposition}

\noindent\textbf{Proof.}
Since $(Y_\cdot, Z_\cdot)$ is a solution of BSDE \eqref{BSDE:4.1}  such that $(e^{\gamma (|Y_t|+\int_0^t\alpha_s\dif s)})_{t\in[0,T]}$ belongs to class $(D)$, according to Lemma \ref{lem-nc} we obtain that there exists a strictly increasing differentiable function $k:\R_+\rightarrow\R_+$ such that $k(0)=\gamma, ~ k(x)\rightarrow+\infty$ when $x\rightarrow+\infty$, and
\begin{align}\label{kk}
\sup_{\tau\in\Sigma_T}\E\Big[K\Big(|Y_\tau|+\int_0^{\tau}\alpha_s\dif s\Big)\Big]<+\infty
\end{align}
with $K(x)=\int_0^x k(u)e^{\gamma u}\dif u,~~x\in\R_+$. Define
\begin{align*}
\Psi(x)=\int_0^x k(u)(e^{\gamma u}-1)\dif u, \ \ x\in\R_+.
\end{align*}
Since $\Psi$ is a convex function on $\R_+$, we know that the dual function of $\Psi$ is
$$\Phi(x)=\int_0^x \Phi'(u)\dif u$$
with $\Phi'$ being the inverse function of $\Psi'$. We consider an $(\F_t)$-progressively measurable process $(q^*_s)_{s\in[0,T]}$ valued in $\R^{1\times d}$ and such that $q_s^*\in\partial_z g_1(s,Z_s)$ for all $s\in[0,T]$. First, we have to show that $\ps$, $\int_0^T |q_s^*|^2\dif s<+\infty$. In fact, since assumptions (A1) and (A2') hold for $g_1$, we can define $f_1$ as in \eqref{fdy1}, and \eqref{fx-nc} holds. It then follows from \eqref{fx-nc} and Young's inequality that $\ps$,
\begin{align*}
|q_s^*|^2&\leq 2\gamma \big(f_1(s,q_s^*)+\alpha_s\big)=2\gamma \big(q_s^*Z_s-g_1(s,Z_s)+\alpha_s\big)\\
&\leq 2\gamma\Big(\frac{1}{4\gamma}|q_s^*|^2+\gamma|Z_s|^2-g_1(s,Z_s)+\alpha_s\Big),\ \ s\in\T.
\end{align*}
Thus,
$$
\int_0^T |q_s^*|^2\dif s\leq 4\gamma\int_0^T \Big(\gamma|Z_s|^2-g_1(s,Z_s)+\alpha_s\Big)\dif s<+\infty.\vspace{0.2cm}
$$
Now, let us show that $\mathcal{E}(q^*)$ is a uniformly integrable martingale. For each $n\geq1$, define the stopping time
$$
\tau_n:=\inf\Big\{t\in[0,T]:\int_0^t |q_s^*|^2\dif s+\int_0^t |Z_s|^2\dif s\geq n\Big\}\wedge T
$$
with the convention $\inf\emptyset=+\infty$, and the probability measure $\Q_n^*$ by
$$
\frac{\dif \Q^*_n}{\dif \mathbb{P}}:=M^*_{\tau_n}\ \ {\rm with}\ \ M^*_t=\exp\Big(\int_0^t q_s^*\dif B_s-\frac{1}{2}\int_0^t |q_s^*|^2\dif s\Big), \ \ t\in\T.
$$
Set $B_t^{q^*}:=B_t-\int_0^t (q^*_s)^{\top}\dif s, ~t\in\T$, then $(B_t^{q^*})_{t\in[0,\tau_n]}$ is a standard $d$-dimensional Brownian motion under the probability $\Q_n^*$ for each $n\geq1$. Now, we verify that $(M^*_{\tau_n})_{n\in\N}$ is uniformly integrable which is sufficient to conclude the desired result. In view of $q_s^*\in\partial_z g_1(s,Z_s)$, thanks to BSDE \eqref{BSDE:4.1}  and Girsanov's theorem, we have
\begin{align}\label{cf2}
\notag Y_0&=Y_{\tau_n}-\int_0^{\tau_n}g(s,Z_s)\dif s+\int_0^{\tau_n} Z_s \cdot \dif B_s\\
\notag &=Y_{\tau_n}+\int_0^{\tau_n}\big(q_s^*Z_s-g(s,Z_s)\big)\dif s+\int_0^{\tau_n} Z_s \cdot \dif B_s^{q^*}\\
\notag &=Y_{\tau_n}+\int_0^{\tau_n}\big(f_1(s,q_s^*)-g_2(s,Z_s)\big)\dif s+\int_0^{\tau_n} Z_s \cdot \dif B_s^{q^*}\\
&=\E^{\Q^*_n}\Big[Y_{\tau_n}+\int_0^{\tau_n}\big(f_1(s,q_s^*)-g_2(s,Z_s)\big)\dif s\Big].
\end{align}
According to assumption (H2) of $g_2$ with $a=0$, we can deduce that
\begin{align*}
|g_2(s,Z_s)|\leq b, ~~~s\in\T,
\end{align*}
which combining the fact that $Y_{\tau_n}\geq-|Y_{\tau_n}|$, \eqref{cf2} and \eqref{fx-nc} gives
\begin{align}\label{bb-nc}
\notag Y_0&\geq\E^{\Q^*_n}\Big[-|Y_{\tau_n}|+\int_0^{\tau_n}\Big(-\alpha_s+\frac{1}{2\gamma}|q_s^*|^2-b\Big)\dif s\Big]\\
&\geq-\E^{\Q^*_n}\Big[|Y_{\tau_n}|+\int_0^{\tau_n}\alpha_s\dif s\Big]+\frac{1}{2\gamma}\E^{\Q^*_n}\Big[\int_0^{\tau_n}|q_s^*|^2\dif s\Big]-bT.
\end{align}
Since $\Psi$ and $\Phi$ are dual convex functions, Young's inequality gives
\begin{align}\label{1xing}
\notag -\E^{\Q^*_n}\Big[|Y_{\tau_n}|+\int_0^{\tau_n}\alpha_s\dif s\Big]&=-\E\Big[\Big(|Y_{\tau_n}|+\int_0^{\tau_n}\alpha_s\dif s\Big)M^*_{\tau_n}\Big]\\
&\geq-\E\Big[\Psi\Big(|Y_{\tau_n}|+\int_0^{\tau_n}\alpha_s\dif s\Big)\Big]-\E[\Phi(M^*_{\tau_n})].
\end{align}
By virtue of \eqref{kk} along with the definitions of $K(\cdot)$ and $\Psi(\cdot)$, we have
\begin{align}\label{2xing}
-\E\Big[\Psi\Big(|Y_{\tau_n}|+\int_0^{\tau_n}\alpha_s\dif s\Big)\Big]\geq -\E\Big[K\Big(|Y_{\tau_n}|+\int_0^{\tau_n}\alpha_s\dif s\Big)\Big]\geq-C,
\end{align}
where $C$ is a constant independent of $n$. Moreover, a simple calculus gives
\begin{align}\label{3xing}
\frac{1}{2\gamma}\E^{\Q^*_n}\Big[\int_0^{\tau_n}|q_s^*|^2\dif s\Big]=\frac{1}{\gamma}\E\left[M^*_{\tau_n}\ln M^*_{\tau_n}\right].
\end{align}
By putting \eqref{1xing}, \eqref{2xing} and \eqref{3xing} into \eqref{bb-nc}, we obtain that
\begin{align*}
\notag Y_0&\geq-C-\E[\Phi(M^*_{\tau_n})]+\frac{1}{\gamma}\E\left[M^*_{\tau_n}\ln M^*_{\tau_n}\right]-bT\\
&=-C-bT+\E[\Lambda(M^*_{\tau_n})],
\end{align*}
where
$$\Lambda(x):=\frac{1}{\gamma}x\ln x-\Phi(x),\ \ x\in\R_+.$$
Hence,
\begin{align}\label{4.28}
\sup_{n\geq1}\E[\Lambda(M^*_{\tau_n})]<+\infty,
\end{align}

Furthermore, we have the following proposition similar to Proposition 2 of \cite{Delbaen2015}. Its proof is given here for readers' convenience.

\begin{proposition}\label{pro-lambda-nc}
The function $\Lambda(x)$ satisfies
$$
\lim_{x\rightarrow+\infty}\frac{\Lambda(x)}{x}=+\infty.
$$
\end{proposition}

\noindent\textbf{Proof.} It is sufficient to show that
$$\Lambda'(x)=\frac{1}{\gamma}\ln x+\frac{1}{\gamma}-\Phi'(x),\ \ x>0$$
is increasing and $\lim_{x\rightarrow+\infty}\Lambda'(x)=+\infty$. First, let us show that $\Psi''(\Phi'(x))\geq\gamma(x+\gamma)$ for all $x\in\R_+$:
$$
\Psi''(x)=k'(x)(e^{\gamma x}-1)+\gamma k(x)e^{\gamma x}\geq \gamma k(x)e^{\gamma x}=\gamma k(x)(e^{\gamma x}-1)+\gamma k(x)\geq\gamma\Psi'(x)+\gamma^2,
$$
so we have
$$
\Psi''(\Phi'(x))\geq\gamma\Psi'(\Phi'(x))+\gamma^2=\gamma(x+\gamma).\vspace{0.3cm}
$$
Thus, in view of $(\Psi'(\Phi'(x)))'=\Psi''(\Phi'(x))\Phi''(x)=1$, we can deduce that
$$
\Phi''(x)\leq\frac{1}{\gamma(x+\gamma)},\ \ x\in\R_+.
$$
Hence, we have
$$
\Lambda''(x)=\frac{1}{\gamma x}-\Phi''(x)\geq\frac{1}{\gamma x}-\frac{1}{\gamma(x+\gamma)}>0,\ \ x>0,\vspace{0.2cm}
$$
which means that $\Lambda'(\cdot)$ is an increasing function on $\R_+$.\vspace{0.2cm}

To conclude we will prove by contradiction that $\Lambda'(\cdot)$ is an unbounded function on $[1,+\infty)$: let us assume that there exists a constant $L>0$ such that $\Lambda'(x)\leq L$ for all $x\in [1,+\infty)$. Then we have
\begin{align*}
x&=\Psi'(\Phi'(x))=k(\Phi'(x))(e^{\gamma\Phi'(x)}-1)=k(\Phi'(x))(e^{\gamma(\frac{1}{\gamma}\ln x+\frac{1}{\gamma}-\Lambda'(x))}-1)\\
&\geq k(\Phi'(x))(e^{\ln x+1}e^{-\gamma L}-1)=k(\Phi'(x))(xe^{1-\gamma L}-1).
\end{align*}
Then, we get for $x$ big enough,
$$
k(\Phi'(x))\leq\frac{x}{xe^{1-\gamma L}-1}\leq C_{\gamma,L},\vspace{0.2cm}
$$
where $C_{\gamma,L}$ is a positive constant depending only on $\gamma$ and $L$. Since $$\lim_{x\rightarrow+\infty}\Phi'(x)=+\infty,\vspace{-0.1cm}$$
the last inequality implies that $k(\cdot)$ is a bounded function, which is a \vspace{0.3cm}contradiction.\hfill\framebox

Finally, according to Lemma \ref{lem-poussin} (the de La Vall$\acute{\rm{e}}$e Poussin lemma) along with \eqref{4.28} and Proposition \ref{pro-lambda-nc}, the conclusion of Proposition \ref{Q} follows \vspace{0.3cm}immediately. \hfill\framebox

Based on Proposition \ref{Q}, we can prove Theorem \ref{ncth2}.
\vspace{0.2cm}

\noindent\textbf{Proof of Theorem \ref{ncth2}.}
Since $g_1$ satisfies (A1), and $g_2$ satisfies (H2) with $a=0$, it is not hard to verify that the generator $g:=g_1+g_2$ of BSDE \eqref{BSDE:4.1}  satisfies (A1). Thus, the existence result in Theorem \ref{ncth2} has been given in (ii) of Proposition \ref{cun}. Now, we are committed to proving the uniqueness. In fact, with Proposition \ref{Q} in hand and in view of $g_1$ satisfying (A1) and (A2'), and $g_2$ satisfying (H1)-(H2) with $a=0$, by an identical argument as that in the proof of (iii) of Theorem \ref{ncth1} we can verify the desired assertion on the uniqueness. The only difference lies in $A=b, q^A_\cdot\equiv 0$ and $B^A_\cdot\equiv B^{q^*}_\cdot$ here. The proof is then \vspace{0.2cm}complete. \hfill\framebox

\begin{remark}
{\rm The novelty of Theorem \ref{ncth2} comes from the fact that the generator is allowed to be a strongly convex function perturbed by a bounded uniformly continuous function, and that the terminal value only needs the sharp exponential moment required for the existence result. This result extends the uniqueness result of [6] where there is no perturbation allowed. It is worth noting that even without the perturbation term, that is, when $a=b=0$ in assumption (H2) (which implies $g_2\equiv 0$), Propositions \ref{Q} and \ref{pro-lambda-nc} still improve the proofs of Propositions 1 and 2 in \cite{Delbaen2015}, \vspace{0.3cm}respectively.}
\end{remark}

%\section*{References}
\setlength{\bibsep}{2pt}

%%-----------------------------
\end{document}